\PassOptionsToPackage{sort}{natbib}
\documentclass[final,3p,times,11pt]{elsarticle}
\usepackage{amsmath,amssymb,amsfonts}
\usepackage{amsthm}
\usepackage[dvipsnames]{xcolor}
\usepackage{dsfont}

\usepackage{graphicx,subcaption,lipsum}
\usepackage{url}
\usepackage{doi}
\usepackage{hyperref}
\usepackage[capitalize]{cleveref}
\usepackage{subcaption}
\usepackage[labelfont=bf]{caption}
\usepackage{float}
\usepackage{orcidlink}
\usepackage{enumitem}
\usepackage{natbib}
\usepackage{tikz}
\usetikzlibrary{decorations.pathreplacing}
\usepackage{pgfplots}
\pgfplotsset{compat=1.18}
\usepackage{algorithm}
\usepackage{algpseudocode}
\usepackage{lineno}

\newtheorem{theorem}{Theorem}[section]
\newtheorem{lemma}[theorem]{Lemma}

\newtheorem{theorem-appendix}{Theorem A.\ignorespaces}%[section]
\newtheorem{lemma-appendix}[theorem-appendix]{Lemma A.\ignorespaces}
\newtheorem{remark-appendix}[theorem-appendix]{Remark A.\ignorespaces}
\renewcommand{\d}[1]{\textnormal{d}#1}

\usepackage{todonotes}
\newboolean{SHOWCOMMENTS}
\setboolean{SHOWCOMMENTS}{true}
\ifthenelse{\boolean{SHOWCOMMENTS}}{
	\newcommand{\tdMK}[1]
	{\todo[color=blue!25,inline]{\footnotesize{\bf Martin:} #1}}
}{
	\newcommand{\tdMK}[1]{}
}

\ifthenelse{\boolean{SHOWCOMMENTS}}{
	\newcommand{\tdAW}[1]
	{\todo[color=teal!25,inline]{\footnotesize{\bf Anna:} #1}}
}{
	\newcommand{\tdAW}[1]{}
}

\journal{}%Mathematics and Computers in Simulation}

\begin{document}

\begin{frontmatter}
	
	
	
	%\title{Applied numerical analysis of higher-order schemes for integro-differential equation-based age-of-infection models with time-dependent contact rates}
	%\title{Development and numerical analysis of Gregory quadrature and finite differences combined higher-order schemes for integro-differential equation-based formulations of age-of-infection models with time-dependent contact rates}
	\title{Higher order schemes for the numerical solution of integro-differential equation-based age-of-infection models with time-dependent contact rates}
	
	\author[label1]{Anna Wendler\hspace{0.5mm}\orcidlink{0000-0002-1816-8907}}
	\ead{anna.wendler@dlr.de}
	\author[label2]{Jan Hasenauer\hspace{0.5mm}\orcidlink{0000-0002-4935-3312}}
	\ead{jan.hasenauer@uni-bonn.de}
	\author[label1,label2]{Martin J. Kühn\hspace{0.5mm}\orcidlink{0000-0002-0906-6984}}
	\ead{martin.kuehn@dlr.de}
	
	\affiliation[label1]{organization={Institute of Software Technology, Department of High-Performance Computing, German Aerospace Center},%Department and Organization, 
		city={Cologne},
		country={Germany}}
	
	\affiliation[label2]{organization={Bonn Center for Mathematical Life Sciences and Life \& Medical Sciences (LIMES) Institute, University of Bonn, Bonn, Germany},%Department and Organization, 
		city={Bonn},
		country={Germany}}
	
	\begin{abstract}
		
		Age-of-infection models, formulated as integro-differential equations (IDEs), provide elegant mathematical descriptions for the spread of infectious diseases. In contrast to ordinary differential equation (ODE)-based models, IDE-based models can accommodate general distributions of transition times between disease stages. 
		Yet, the accurate numerical assessment of the solution of IDE models remains challenging.
		
		Here, we extend a recently proposed fixed-point method based on Gregory quadrature. We generalize the formulation to time-dependent contact rates whose changes at discrete time points are represented by continuously differentiable sigmoidal interpolants, and suggest additional finite difference approximations to obtain higher-order schemes for susceptible, infected and recovered subpopulations. The resulting numerical schemes achieve convergence orders of two, three and four for all three subpopulations, enabling accurate and flexible simulations of IDE-based epidemic models with general sojourn time distributions and changing contact patterns.
	\end{abstract}
	
	\begin{highlights}
		\item Time-dependent contact rates for age-of-infection models 
		\item Higher order numerical schemes including infected and recovered populations
		\item Combination of Gregory quadrature rules with finite difference approximations
		\item Numerical validation of higher order convergence for discretized model
		\item Assessment of extended model with respect to its impact on epidemic outcomes
	\end{highlights}
	
	
\end{frontmatter}
% \linenumbers

\section{Introduction}\label{sec:introduction}

Ordinary differential equation (ODE) models are among the most popular tools for describing and analyzing infectious disease dynamics in a population~\cite{hoang2025differential}. These models allow for the computation of reproduction numbers and (endemic) equilibria; however, they implicitly assume exponentially distributed sojourn times in the different disease states, an assumption that is often unrealistic~\cite{feng_epidemiological_2007}. More complex sojourn-time distributions can be represented within the ODE framework using the Linear Chain Trick (LCT), which yields Erlang-distributed transition times~\cite{macdonald_time_1978}, or its generalizations~\cite{hurtado_generalizations_2019,hurtado_building_2021,hurtado_finding_2025}. Although these approaches provide greater modeling flexibility and can reveal dynamics absent from simpler ODE models~\cite{ploetzke_lct_2025}, they require the introduction of auxiliary states. This increases the dimensionality of the model and can lead to intricate systems, particularly when complex sojourn-time distributions are considered. 

A mathematically direct and more elegant alternative to expanded ODE models is provided by integro-differential equation-based (IDE) models which accommodate general sojourn-time distributions without introducing auxiliary states. IDE-based epidemic models were already considered by Kermack and McKendrick~\cite{kermackContribution1927} and by Hethcote and Tudor~\cite{hethcote_integral_1980}. These models have also commonly been referred to as age-of-infection models~\cite{brauer_age_2009,brauer_mathematical_2019,barrilReproductionNumberAge2021,demongeotKermackMcKendrickModel2023,pezzellaInfectiousDisease2024}.

A substantial body of analytical work has established the theoretical foundations for SIR-type age-of-infection models. These contributions include analyses of the classical Kermack--McKendrick model~\cite{brauerKermackMcKendrickEpidemic2005}, final-size relations~\cite{brauerAgeinfectionFinalSize2008}, reproduction numbers and final sizes under heterogeneous mixing~\cite{brauer_age_2009}, endemic equilibria~\cite{breda_formulation_2012}, and the dynamics of disease-specific models~\cite{brauerDynamicsAgeinfectionCholera2013}. 
More recently, the COVID-19 pandemic motivated the development of more advanced models that distinguish between symptomatic and asymptomatic infections~\cite{baiAgeofinfectionModelBoth2023,wendler_nonstandard_2026}. Further work also demonstrated that the long-term dynamics can depend on the distributional shape of the infectious period~\cite{scarabel_bistability_2026}. Despite this extensive foundation, important questions concerning the accurate and efficient numerical solution of nonlinear IDE-based epidemic models remain open. The lack of established numerical methods represents a practical barrier to the broader use of these models, particularly when compared with the wide range of mature integration schemes available for ODE models.

Numerical methods for IDE-based epidemic models are a field of active research, and the absence of user-friendly computational tools remains a recognized limitation~\cite{Diekmann2025}. An early contribution was the nonstandard scheme developed by Lubuma and Terefe for an SIS model~\cite{lubuma_nonstandard_2015}. Building on this, Messina, Pezzella, and Vecchio introduced a nonstandard finite difference scheme that applies modified left- and right-rectangle rules to different terms of the integral equation~\cite{messina_non-standard_2022,pezzellaInfectiousDisease2024}. This construction resolves the nonlinear coupling between the susceptible population and the force of infection while preserving positivity, monotonicity, and boundedness independently of the step size $h$. The scheme was subsequently applied to models incorporating heterogeneous population mixing~\cite{messina_longtime_behavior_2024} and behavioral changes in contact patterns~\cite{buonomo_behavioral_changes_2024}, and it underlies the EPITIME computational framework for behavioral epidemic modeling~\cite{buonomo_epitime_2026}. We further extended the approach to a SECIR-type IDE model that includes disease-related deaths and established preservation properties for its eight compartments and the flows between them~\cite{wendler_nonstandard_2026}. Despite these advantages, the nonstandard scheme is only first-order convergent.

Enabling efficient, high-accuracy simulations of complex models demands higher-order methods. Promising approaches have been given by~\cite{messina_positive_2022,pezzellaInfectiousDisease2024}, where the authors combine a fixed-point iteration with Gregory's quadrature method for the differential equation of the susceptible population. The Gregory quadrature weights used for methods of order $n\in\mathbb{N}$ are obtained when adding further correction terms to the composite trapezoidal rule; see also~\cite{brass1977quadraturverfahren}.

In this paper, we consider a generalized version of the continuous IDE-based model of~\cite{messina_positive_2022,pezzellaInfectiousDisease2024}, allowing for time-dependent contact rates and the isolation of infected individuals as in~\cite{wendler_nonstandard_2026} for which we extend the numerical scheme accordingly. More explicitly, we examine the convergence for scenarios with changing contact rates, e.g., due to nonpharmaceutical interventions and apply $k$-times continuously differentiable interpolation functions to smooth the transition step. Using the fixed-point iteration of~\cite{messina_positive_2022,pezzellaInfectiousDisease2024} and higher-order finite differences, we develop combined higher-order methods, which are not restricted to but numerically validated for orders $n\in\{2,3,4\}$, for the differential equations describing the susceptible, infected, and recovered population.

The paper is structured as follows: In~\cref{sec:continuous_model}, we introduce the generalized continuous IDE formulation. In~\cref{sec:discretization}, we present the discretization of the continuous formulation through our suggested numerical scheme. In~\cref{sec:numerical_results}, we provide numerical results for different aspects of a ground truth solution and the suggested numerical scheme before drawing a conclusion in~\cref{sec:conclusion}.

\section{An age-of-infection model through a continuous formulation with integro-differential equations}\label{sec:continuous_model}

\subsection{Preliminaries and general description of dynamics}

We consider a continuous SIR-type epidemic model that subdivides the population into the following three, time-dependent compartments:
\begin{itemize}
	\item \textit{Susceptible} individuals $S:\mathbb{R}\rightarrow \mathbb{R}_{\geq 0}$ who have not been infected with the disease yet but are infected upon contact with a nonisolated infected individual with probability $\rho$ (that may depend on the time since infection);
	\item \textit{Infected} individuals $I:\mathbb{R}\rightarrow \mathbb{R}_{\geq 0}$ who are infected and infectious;
	\item \textit{Recovered} individuals $R:\mathbb{R}\rightarrow \mathbb{R}_{\geq 0}$ who are immune through prior infection and cannot be infected again.
\end{itemize}

We define the contact rate $\phi:\mathbb{R}\rightarrow\mathbb{R}_{\geq 0}$, the transmission probability $\rho:\mathbb{R}_{\geq 0}\rightarrow[0,1]$, and the fraction of infected individuals that is not isolated as $\xi:\mathbb{R}_{\geq 0}\rightarrow[0,1]$. While $\phi$ is evaluated along with the compartments, the functions $\rho$ and $\xi$ depend on the time since infection. Correspondingly, $\gamma_I^R:\mathbb{R}_{\geq0}\rightarrow[0,1]$ denotes the survival function of the sojourn time distribution in the infected compartment. It is monotonically decreasing, satisfies $\gamma_I^R(0)=1$, and is assumed to be continuously differentiable on $(0,\infty)$, with
\begin{align*}
	\int_0^\infty\gamma_I^R(s)\,\d s<\infty.
\end{align*}
Thus, the expected sojourn time is finite and individuals eventually leave the infected compartment. These infection age-dependent functions also govern the future evolution of individuals infected before time zero.

With $\lim_{t\rightarrow-\infty} S(t)=N$, we define the change in the susceptible population by
\begin{align}\label{eq:S_cont_infty}
	S'(t)&=\frac{\phi(t)}{N}S(t)
	\int_0^\infty\rho(s)\,\xi(s)\,\gamma_I^R(s)\,S'(t-s)\,\d s.
\end{align}
The corresponding cohort expressions for the number of individuals in infected and recovered state are
\begin{align}\label{eq:I_R_cont_infty}
	\begin{aligned}
		I(t)&=-\int_{-\infty}^t\gamma_I^R(t-s)\,S'(s)\,\d s,\\
		R(t)&=-\int_{-\infty}^t\bigl(1-\gamma_I^R(t-s)\bigr)\,S'(s)\,\d s.
	\end{aligned}
\end{align}
%These expressions assume that all infections are represented in the historical incidence and that\linebreak $\lim_{s\to-\infty}S(s)=N$; pre-existing immune populations or separately introduced seed cohorts require additional terms. Prescribed initial values for the compartments at time zero must agree with these cohort expressions. Yet, this consistency requirement does not require the historical incidence to satisfy the transmission equation before time zero.

For a shorter notation, we follow~\cite{messina_non-standard_2022} and set
\begin{align*}
	A(s)=\rho(s)\,\xi(s)\,\gamma_I^R(s),\qquad
	\beta(t)=\frac{\phi(t)}{N},
\end{align*}
to also use
\begin{align}\label{eq:S_cont_infty_A}
	S'(t)&=\beta(t)\,S(t)\int_{-\infty}^t A(t-s)\,S'(s)\,\d s.
\end{align}

\subsection{Finite-history approximation}\label{sec:finite_history}

We assume that the total population $N\in\mathbb{N}$ is constant and choose a simulation start $\widetilde{T}\leq 0$, which does not necessarily coincide with the start of the epidemic. To obtain a finite-history formulation, we approximate the infection-age distribution at $\widetilde{T}$ by assigning infection age zero to all individuals infected at that time. The resulting approximation error generally depends on the choice of $\widetilde{T}$. For the initial values, we require $S(\widetilde{T})+I(\widetilde{T})+R(\widetilde{T})=N$. The case $I(\widetilde{T})=0$ is trivial, with the system remaining in the disease-free equilibrium $(S(\widetilde{T}),0,R(\widetilde{T}))$. Commonly, the simulation start is chosen as $\widetilde{T}=0$ together with $R(\widetilde{T})=0$; see~\cite{brauerKermackMcKendrickEpidemic2005,pezzellaInfectiousDisease2024}.

Mathematically, we proceed as follows and split the integral in~\eqref{eq:S_cont_infty_A} as
\begin{align}\label{eq:history_split}
	S'(t)=&\,\beta(t)\,S(t)\left(
	\int_{\widetilde{T}}^t A(t-s)\,S'(s)\,\d s+B_{\widetilde{T}}(t)\right),\\
	B_{\widetilde{T}}(t):=&\int_{-\infty}^{\widetilde{T}}A(t-s)\,S'(s)\,\d s.
\end{align}
The function $B_{\widetilde{T}}$ contains the contribution from infections predating $\widetilde{T}$. To approximate this contribution, we proceed as~\cite{brauerAgeinfectionFinalSize2008,messina_positive_2022} and assume that all infected at time $\widetilde{T}$ have infection age zero. However, in order to allow simulation starts after the occurrence of the first infection(s) of the epidemic, we need to take into account potential recoveries. Thus, we have
\begin{align}\label{eq:assump_initial_conditions}
	\begin{aligned}
		B_{\widetilde{T}}(t) &\approx A(t-\widetilde{T})\left(\int_{t-\widetilde{T}}^\infty \,S'(t-s)\,\d s + R(\widetilde{T})\right)\\
		&=\left(S(\widetilde{T})-N+R(\widetilde{T})\right)\,A(t-\widetilde{T})\\
		&=-I(\widetilde{T})\,A(t-\widetilde{T}),
	\end{aligned}
\end{align}
This representation is exact for a cohort that is newly infected at $\widetilde{T}$. For a population with a nontrivial infection-age distribution at $\widetilde{T}$, it is generally an approximation to the older history. With $\widetilde{T}=0$ and $R(\widetilde{T})=0$, it yields the initialization considered in~\cite[Eq.~2]{messina_positive_2022}.

The error describing the difference between the full and approximated historical contributions is
\begin{align}\label{eq:history_error_Sprime}
	E_{S',\widetilde{T}}(t)
	:=B_{\widetilde{T}}(t)+I(\widetilde{T})\,A(t-\widetilde{T}).
\end{align}
For the sake of a simplified presentation, we do not distinguish the notation for the functions $S, I, R$ defined on $\mathbb{R}$ through~\eqref{eq:S_cont_infty}~and~\eqref{eq:I_R_cont_infty} and their corresponding counterparts defined through the following system using Assumption~\eqref{eq:assump_initial_conditions}, i.e., on the left-closed interval with finite history $[\widetilde{T},t)$. The finite-history model with the derivative for $S(t)$ obtained using~\eqref{eq:assump_initial_conditions} is given by
\begin{align}\label{eq:Sprime_finite}
	S'(t)&=\beta(t)\,S(t) \left(
	\int_{\widetilde{T}}^t A(t-s)\,S'(s)\,\d s-I(\widetilde{T})\,A(t-\widetilde{T})\right)
\end{align}
and
\begin{align}\label{eq:I_R_cont_finite}
	\begin{aligned}
		I(t)
		&=\gamma_I^R(t-\widetilde{T})\,I(\widetilde{T})
		-\int_{\widetilde{T}}^t
		\gamma_I^R(t-s)\,S'(s)\,\mathrm{d}s,\\
		R(t)
		&=R(\widetilde{T})
		+\bigl(1-\gamma_I^R(t-\widetilde{T})\bigr)\,I(\widetilde{T})
		-\int_{\widetilde{T}}^t
		\bigl(1-\gamma_I^R(t-s)\bigr)\,S'(s)\,\mathrm{d}s\\
		&=N-S(t)-I(t).
	\end{aligned}
\end{align}
Here, the first term in $I(t)$ describes survivors of the cohort assigned infection age zero at $\widetilde{T}$; the integral describes survivors of subsequent infections. Analogously to~\eqref{eq:history_error_Sprime}, we can define
\begin{align}\label{eq:error_IR}
	\begin{aligned}
		E_{I,\widetilde{T}}(t)&:=\int_{-\infty}^{\widetilde{T}}\gamma_I^R(t-s)\,S'(s)\,\d s + \gamma_I^R(t-\widetilde{T})\,I(\widetilde{T})\\
		\text{and}\quad E_{R,\widetilde{T}}(t)&:=\int_{-\infty}^{\widetilde{T}}(1-\gamma_I^R(t-s))\,S'(s)\,\d s + R(\widetilde{T}) + (1-\gamma_I^R(t-\widetilde{T}))\,I(\widetilde{T})\\
	\end{aligned}
\end{align}

\subsection{Particular results for exponential sojourn time distribution}\label{sec:exp}

While our motivation for IDE formulations is explicitly based on nonexponential sojourn times, an important special case is exponential survival, $\gamma_I^R(s)=e^{-s/T_I}$, with $T_I$ the expected sojourn time in the infected compartment and where, additionally, $\rho(s)=\rho$ and $\xi(s)=\xi$ are assumed to be infection age-independent. We now formulate two simple lemmas which will be used in the numerical results.

\begin{lemma}\label{lemma:error}
	Consider the model of~\eqref{eq:S_cont_infty} and~\eqref{eq:I_R_cont_infty} and assume $\gamma_I^R(s)=e^{-s/T_I}$ as well as $\rho(s)=\rho$ and $\xi(s)=\xi$. Then, no error is introduced when reducing the model as in~\eqref{eq:S_cont_infty} and~\eqref{eq:I_R_cont_infty} to its finite history approximation as in~\eqref{eq:Sprime_finite} and~\eqref{eq:I_R_cont_finite}.
\end{lemma}
\begin{proof}
	We only need to show that the error terms~\eqref{eq:history_error_Sprime} and~\eqref{eq:error_IR} vanish. In the case of exponentially distributed sojourn time distributions, the two identities
	\begin{align*} 
		\gamma_I^R(u+v)&=\gamma_I^R(u)\,\gamma_I^R(v)\qquad \text{and}\qquad
		\gamma_I^R(-u)= \frac{1}{\gamma_I^R(u)}
	\end{align*}
	hold. With the first identity, we obtain from~\eqref{eq:I_R_cont_infty} that
	\begin{align*}
		I(\widetilde{T})&=-\int_{-\infty}^{\widetilde{T}}\gamma_I^R(\widetilde{T}-s)\,S'(s)\,\d s\\
		&= -\int_{-\infty}^{\widetilde{T}}\gamma_I^R(\widetilde{T}-t)\,\gamma_I^R(t-s)\,S'(s)\,\d s\\
		&= -\gamma_I^R(\widetilde{T}-t)\,\int_{-\infty}^{\widetilde{T}}\gamma_I^R(t-s)\,S'(s)\,\d s.
	\end{align*}
	Multiplying with constant $-\rho\,\xi$ and using the second identity on $\gamma_I^R(t-\widetilde{T})$ yields
	\begin{align*}
		-\underset{=A(t-\widetilde{T})}{\underbrace{\rho\,\xi\,\gamma_I^R(t-\widetilde{T})}}\,I(\widetilde{T}) = \int_{-\infty}^{\widetilde{T}}\underset{=A(t-s)}{\underbrace{\rho\,\xi\,\gamma_I^R(t-s)}}\,S'(s)\,\d s = B_{\widetilde{T}}(t),
	\end{align*}
	which implies $E_{S',\widetilde{T}}(t)=0$. With similar arguments as above, we also obtain $E_{I,\widetilde{T}}(t)=0$ and
	\begin{align*}
		E_{R,\widetilde{T}}(t)&=E_{I,\widetilde{T}}(t)+\int_{-\infty}^{\widetilde{T}}S'(s)\,\d s + R(\widetilde{T}) + I(\widetilde{T})
		=S(\widetilde{T})-\lim_{t\rightarrow-\infty} S(t)+ R(\widetilde{T}) + I(\widetilde{T}),
	\end{align*}
	which vanishes with our assumption on the initial population $S(\widetilde{T})+I(\widetilde{T})+R(\widetilde{T})=N$ and $\lim_{t\rightarrow-\infty} S(t)=N$.
\end{proof}
The consequence of~\cref{lemma:error} is that there is no difference between the outcome of an IDE model with full or approximated historical contributions when using exponential sojourn time distributions.

\begin{lemma}\label{lemma:ode}
	Consider the model of~\eqref{eq:Sprime_finite} and~\eqref{eq:I_R_cont_finite} and assume $\gamma_I^R(s)=e^{-s/T_I}$ as well as $\rho(s)=\rho$ and $\xi(s)=\xi$. Then,~\eqref{eq:Sprime_finite} and~\eqref{eq:I_R_cont_finite} can be simplified to the ODE formulation
	\begin{align}\label{eq:SIR_ODE}
		\begin{aligned}
			S'(t)&=-\frac{\phi(t)}{N}S(t)\,
			\rho\,\xi\,I(t),\\
			I'(t)&=\frac{\phi(t)}{N}S(t)
			\rho\,\xi\,I(t)-\frac{1}{T_I}I(t),\\
			R'(t)&=\frac{1}{T_I}I(t).
		\end{aligned}
	\end{align}
\end{lemma}
\begin{proof}
	The proof is obtained by straightforward computation using the Leibniz integral rule. 
\end{proof}

\subsection{Analytical reformulation of the expression for susceptible individuals}

\paragraph{A direct expression for $S$}
In order to obtain an expression that is not given for the derivative but the function $S(t)$ itself, we use an exponential operator. Assuming $\beta(t)$ to be continuously differentiable and adopting the finite-history approximation, we can use a variable transformation and the Leibniz integral rule to rewrite~\eqref{eq:Sprime_finite} as
\begin{align}\label{eq:start_reformulation}
	\begin{aligned}
		S'(t) &= \beta(t)\,S(t) \left(\frac{\d}{\d t}\int_0^{t-\widetilde{T}} A(s)\,S(t-s) \, \d s + \left(R(\widetilde{T})-N\right) \,A(t-\widetilde{T})\right).
	\end{aligned}
\end{align}

When the contact rate is time-dependent, the model introduces additional terms compared to the constant contact rate considered in~\cite{messina_positive_2022}. To handle these terms, we begin with the following observation 
\begin{align}\label{eq:product_deriv}
	\frac{\d}{\d t}\,\left(\beta(t)\int_0^{t-\widetilde{T}} A(s)\,S(t-s) \, \d s\right) &= \beta'(t)\int_0^{t-\widetilde{T}} A(s)\,S(t-s) \, \d s + \beta(t)\,\frac{\d}{\d t}\int_0^{t-\widetilde{T}} A(s)\,S(t-s) \, \d s.
\end{align}

Analogously to~\cite{messina_positive_2022}, we start from~\eqref{eq:start_reformulation} and divide by $S(t)$, for which we require the initial value to be positive, i.e., $S(\widetilde{T})>0$. Integrating from $\widetilde{T}$ to $t$, using~\eqref{eq:product_deriv} and multiple variable transformations of the form $x=t-s$ and $x-\widetilde{T}=t-s$ yields
\begin{align}\label{eq:SprimebyS}
	\begin{aligned}
		\int_{\widetilde{T}}^t \frac{S'(x)}{S(x)}\, \d x 
		&= \int_{\widetilde{T}}^t \beta(x) \frac{\d}{\d x}\int_0^{x-\widetilde{T}} A(s)\,S(x-s) \, \d s \, \d x + \int_{\widetilde{T}}^t \beta(x)\,A(x-\widetilde{T})\,\left(R(\widetilde{T})-N\right) \, \d x \\
		&\hspace{-1.5cm}\overset{\eqref{eq:product_deriv}}{=}  \int_{\widetilde{T}}^t \frac{\d}{\d x}\left(\beta(x)\int_0^{x-\widetilde{T}} A(s)\,S(x-s) \, \d s\right)\ \d x -\int_{\widetilde{T}}^t \beta'(x)\int_0^{x-\widetilde{T}} A(s)\,S(x-s) \, \d s \, \d x\\
		&\hspace{0.5cm}+ \int_{\widetilde{T}}^t \beta(x)\,A(x-\widetilde{T})\,\left(R(\widetilde{T})-N\right) \, \d x \\
		&\hspace{-1.5cm}= \beta(t)\int_0^{t-\widetilde{T}} A(s)\,S(t-s) \, \d s %- \beta(\widetilde{T})\int_0^{\widetilde{T}-\widetilde{T}} A(s)\,S(\widetilde{T}-s) \, \d s \\
		- \int_{\widetilde{T}}^t \beta'(x)\int_0^{x-\widetilde{T}} A(s)\,S(x-s) \, \d s \, \d x + \int_{\widetilde{T}}^t \beta(x)\,A(x-\widetilde{T})\left(R(\widetilde{T})-N\right) \, \d x\\
		&\hspace{-1.5cm}= \beta(t)\int_{\widetilde{T}}^t A(t-s)\, S(s) \, \d s  - \int_{\widetilde{T}}^t \beta'(x)\int_{\widetilde{T}}^x A(x-s)\,S(s) \, \d s \, \d x + \int_{\widetilde{T}}^t \beta(t-s+\widetilde{T})A(t-s)\,\left(R(\widetilde{T})-N\right) \, \d s\\
		&\hspace{-1.5cm}= \int_{\widetilde{T}}^t A(t-s) \left(\beta(t)\,S(s) + \beta(t-s+\widetilde{T})\,\left(R(\widetilde{T})-N\right)\right) \, \d s - \int_{\widetilde{T}}^t \beta'(x)\int_{\widetilde{T}}^x A(x-s)\,S(s) \, \d s \, \d x.
	\end{aligned}
\end{align}

By substituting the derivative within the integral and applying the substitution rule, we obtain the solution for the integral on the left hand side of~\eqref{eq:SprimebyS}. Separating $S(t)$ 
% and re-inserting our expressions for the functions $\beta$ and $A$ 
leads to 
\begin{align}\label{eq:final_S_cont_finite}
	\begin{aligned}
		S(t) &= S(\widetilde{T})\exp\left(\int_{\widetilde{T}}^t A(t-s)\Big(\beta(t)\,S(s)+ \beta(t-s+\widetilde{T})\,\left(R(\widetilde{T})-N\right)\Big)\, \d s - \int_{\widetilde{T}}^t \beta'(x)\int_{\widetilde{T}}^x A(x-s)\,S(s) \, \d s\, \d x\right).\\    
		%&=S(\widetilde{T})\exp\Bigg(\int_{\widetilde{T}}^t \rho(t-s)\,\xi(t-s)\,\gamma_I^R(t-s) \Bigg(\frac{\phi(t)}{N}\,S(s) + \phi(t-s+\widetilde{T})\left(\frac{R(\widetilde{T})}{N}-1\right)\Bigg) \, \d s  \\
		% &\hspace{2cm}- \int_{\widetilde{T}}^t \frac{\phi'(x)}{N}\int_{\widetilde{T}}^x \rho(x-s)\,\xi(x-s)\,\gamma_I^R(x-s)\,S(s) \, \d s \, \d x\Bigg).
	\end{aligned}
\end{align}

\paragraph{An explicit expression for \(S'\)}
For the computation of $I(t)$ and $R(t)$, we require the derivative $S'(t)$. Assuming $A(t)$ to be continuously differentiable, we provide an explicit expression based on the availability of $S(t)$.

For convenience, we define
\begin{align*}
	C:=R(\widetilde{T})-N,
\end{align*}
and denote the exponent in~\eqref{eq:final_S_cont_finite} by
\begin{align}
	H(t)
	&:=
	\int_{\widetilde{T}}^t
	A(t-s)\,
	\Bigl(
	\beta(t)\,S(s)
	+\beta(t-s+\widetilde{T})\,C
	\Bigr)
	\,\d s
	-
	\int_{\widetilde{T}}^t
	\beta'(x)
	\int_{\widetilde{T}}^x A(x-s)\,S(s)\,\d s\,\d x.
	\label{eq:def_H_for_Sprime}
\end{align}
Thus,
\[
S(t)=S(\widetilde{T})\exp\bigl(H(t)\bigr)
\]
and consequently
\begin{align}
	S'(t)=S(t)\,H'(t).
	\label{eq:Sprime_from_H}
\end{align}

Applying the Leibniz integral rule
to~\eqref{eq:def_H_for_Sprime} yields
\begin{align*}%\label{eq:Hprime_intermediate}
	\begin{aligned}
		H'(t)
		&=
		A(0)\,
		\Bigl(
		\beta(t)\,S(t)+\beta(\widetilde{T})\,C
		\Bigr)
		+
		\int_{\widetilde{T}}^t
		A'(t-s)\,
		\Bigl(
		\beta(t)\,S(s)
		+\beta(t-s+\widetilde{T})\,C
		\Bigr)
		\,\d s\\
		&\hspace*{1cm}
		+
		\int_{\widetilde{T}}^t
		A(t-s)\,
		\Bigl(
		\beta'(t)\,S(s)
		+\beta'(t-s+\widetilde{T})\,C
		\Bigr)
		\,\d s
		-
		\beta'(t)
		\int_{\widetilde{T}}^t A(t-s)\,S(s)\,\d s.
	\end{aligned}
\end{align*}
The terms containing
\[
\beta'(t)\int_{\widetilde{T}}^t A(t-s)\,S(s)\,\d s
\]
cancel. Therefore,
\begin{align*}
	%\label{eq:Hprime_before_product}
	\begin{aligned}
		H'(t)
		&=
		A(0)\,\beta(t)S(t)
		+
		\beta(t)\int_{\widetilde{T}}^t A'(t-s)\,S(s)\,\d s\\
		&\quad
		+
		C\Bigl(
		A(0)\,\beta(\widetilde{T})
		+
		\int_{\widetilde{T}}^t
		\left(
		A'(t-s)\,\beta(t-s+\widetilde{T})
		+
		A(t-s)\,\beta'(t-s+\widetilde{T})
		\right)
		\,\d s
		\Bigr).
	\end{aligned}
\end{align*}
Using the substitution \(u=t-s\), the remaining integral satisfies
\begin{align*} %\label{eq:product_integral_Sprime}
	\begin{aligned}
		\int_{\widetilde{T}}^t
		\left(
		A'(t-s)\,\beta(t-s+\widetilde{T})
		+
		A(t-s)\,\beta'(t-s+\widetilde{T})
		\right)
		\,\d s
		& =
		\int_0^{t-\widetilde{T}}
		\frac{\d}{\d u}
		\Bigl(A(u)\,\beta(u+\widetilde{T})\Bigr)
		\,\d u\\
		&
		=
		A(t-\widetilde{T})\,\beta(t)-A(0)\,\beta(\widetilde{T}).
	\end{aligned}
\end{align*}
The term \(A(0)\,\beta(\widetilde{T})\,C\) consequently cancels as well and
we obtain
\begin{equation}
	H'(t)
	=
	\beta(t)\,
	\Bigl(
	A(0)\,S(t)
	+
	\int_{\widetilde{T}}^t A'(t-s)\,S(s)\,\d s
	+
	C\, A(t-\widetilde{T})
	\Bigr).
	\label{eq:Hprime_final}
\end{equation}
Combining~\eqref{eq:Sprime_from_H}
and~\eqref{eq:Hprime_final} gives the explicit expression
\begin{equation}
	S'(t)
	=
	\beta(t)\,S(t)\,
	\Bigl(
	A(0)\,S(t)
	+
	\int_{\widetilde{T}}^t A'(t-s)\,S(s)\,\d s
	+
	C\,A(t-\widetilde{T})
	\Bigr),
	\label{eq:Sprime_explicit}
\end{equation}
where the right-hand side does not depend on $S'(t)$. 

We can easily show that the above formulation is consistent with the original formulation. Assuming that $S(t)$ is sufficiently smooth, integration by parts gives
\begin{align}
	\int_{\widetilde{T}}^t A'(t-s)\,S(s)\,\d s
	&=
	-A(0)\,S(t)
	+A(t-\widetilde{T})\,S(\widetilde{T})
	+\int_{\widetilde{T}}^t A(t-s)\,S'(s)\,\d s.
	\label{eq:integration_by_parts_Sprime}
\end{align}
Using $S(\widetilde{T})+I(\widetilde{T})+R(\widetilde{T})=N$, the expression in parentheses in~\eqref{eq:Sprime_explicit} becomes
\begin{align*}
	&A(t-\widetilde{T})\,\bigl(S(\widetilde{T})+R(\widetilde{T})-N\bigr)
	+
	\int_{\widetilde{T}}^t A(t-s)\,S'(s)\,\d s
	=    -I(\widetilde{T})A(t-\widetilde{T})
	+
	\int_{\widetilde{T}}^t A(t-s)S'(s)\,\d s.
\end{align*}

\subsection{Modeling the contact rate}% DONE/MERGED WITH SUGGESTION

An essential contribution of this work is the systematic accommodation of time-dependent contact rates, which are critical for modeling nonpharmaceutical interventions~\cite{prem_effect_2020} but often disregarded in analytical IDE studies. While our theoretical and numerical framework naturally supports any sufficiently differentiable contact rate function $\phi \in \mathcal{C}^k(\mathbb{R})$, realistic intervention scenarios such as sudden venue closures often involve abrupt, non-smooth changes. Treating these as raw step functions would degrade the convergence order of our high-accuracy numerical schemes; consequently, explicit smoothing must be implemented.

\begin{figure}
	\def\tOne{2}
	\begin{tikzpicture}
		% ------------------------------------------------------------
		%  Cosine transition
		% ------------------------------------------------------------
		\begin{axis}[
			name=polyplot,
			width=7cm,
			height=6cm,
			xmin=0,
			xmax=5,
			ymin=2,
			ymax=8,
			axis lines=left,
			xlabel={$t$},
			ylabel={$\phi(t)$},
			title={Cosine-based interpolation},
			xtick={\tOne,{(\tOne+1)}},
			xticklabels={$t_D - \delta$,$t_D$},
			ytick={3,7},
			yticklabels={$(1-r)\,\phi_0$,$\phi_0$},
			grid=major,
			samples=200,
			domain=0:5,
			]
			
			% Constant value before t_D
			\addplot[
			very thick,
			domain=0:\tOne
			]
			{7};
			
			% Cosine transition:
			% 0.5(4/7 * 7) cos(pi(x-t_D)) + 0.5(14-4)
			% = 2 cos(pi(x-t_D)) + 5.
			%
			% PGFPlots uses degrees for trigonometric functions,
			% so cos(deg(pi*(x-t_D))) is used.
			\addplot[
			very thick,
			domain=\tOne:{\tOne+1},
			]
			{2*cos(deg(pi*(x-\tOne))) + 5};
			
			% Constant value after t_D+1
			\addplot[
			very thick,
			domain={\tOne+1}:5
			]
			{3};
			
			\draw[dashed]
			(axis cs:\tOne,2)
			--
			(axis cs:\tOne,8);
			
			\draw[dashed]
			(axis cs:{\tOne+1},2)
			--
			(axis cs:{\tOne+1},8);
			
		\end{axis}
		
		% ------------------------------------------------------------
		% Ninth-degree smooth transition
		% ------------------------------------------------------------
		\begin{axis}[
			at={(polyplot.east)},
			anchor=west,
			xshift=3cm,
			width=7cm,
			height=6cm,    
			xmin=0,
			xmax=5,
			ymin=2,
			ymax=8,
			axis lines=left,
			xlabel={$t$},
			ylabel={$\phi(t)$},
			title={Polynomial interpolation},
			xtick={\tOne,{(\tOne+1)}},
			xticklabels={$t_D - \delta$,$t_D$},
			ytick={3,7},
			yticklabels={$(1-r)\,\phi_0$,$\phi_0$},
			grid=major,
			samples=200,
			domain=0:5,
			]
			
			% Constant value before t_D
			\addplot[
			very thick,
			domain=0:\tOne
			]
			{7};
			
			% Smooth transition on [t_D,t_D+1].
			% Since x_right-x_left=1, the normalized variable is t=x-t_D.
			\addplot[
			very thick,
			domain=\tOne:{\tOne+1},
			]
			{
				7
				- 4*(
				126*(x-\tOne)^5
				- 420*(x-\tOne)^6
				+ 540*(x-\tOne)^7
				- 315*(x-\tOne)^8
				+ 70*(x-\tOne)^9
				)
			};
			
			% Constant value after t_D+1
			\addplot[
			very thick,
			domain={\tOne+1}:5
			]
			{3};
			
			\draw[dashed]
			(axis cs:\tOne,2)
			--
			(axis cs:\tOne,8);
			
			\draw[dashed]
			(axis cs:{\tOne+1},2)
			--
			(axis cs:{\tOne+1},8);
			
		\end{axis}
	\end{tikzpicture}
	\caption{\textbf{Two sigmoidal interpolation functions to smooth a stepwise defined function between $t_D - \delta$ and $t_D$.} Left: The \textit{cosine-based interpolation} of~\eqref{eq:smooth_cos} ensures $\mathcal{C}^1$ smoothness. Right: The \textit{polynomial interpolation} of~\eqref{eq:smooth_pol} ensures $\mathcal{C}^4$ smoothness.}
	\label{fig:smoother}
\end{figure}
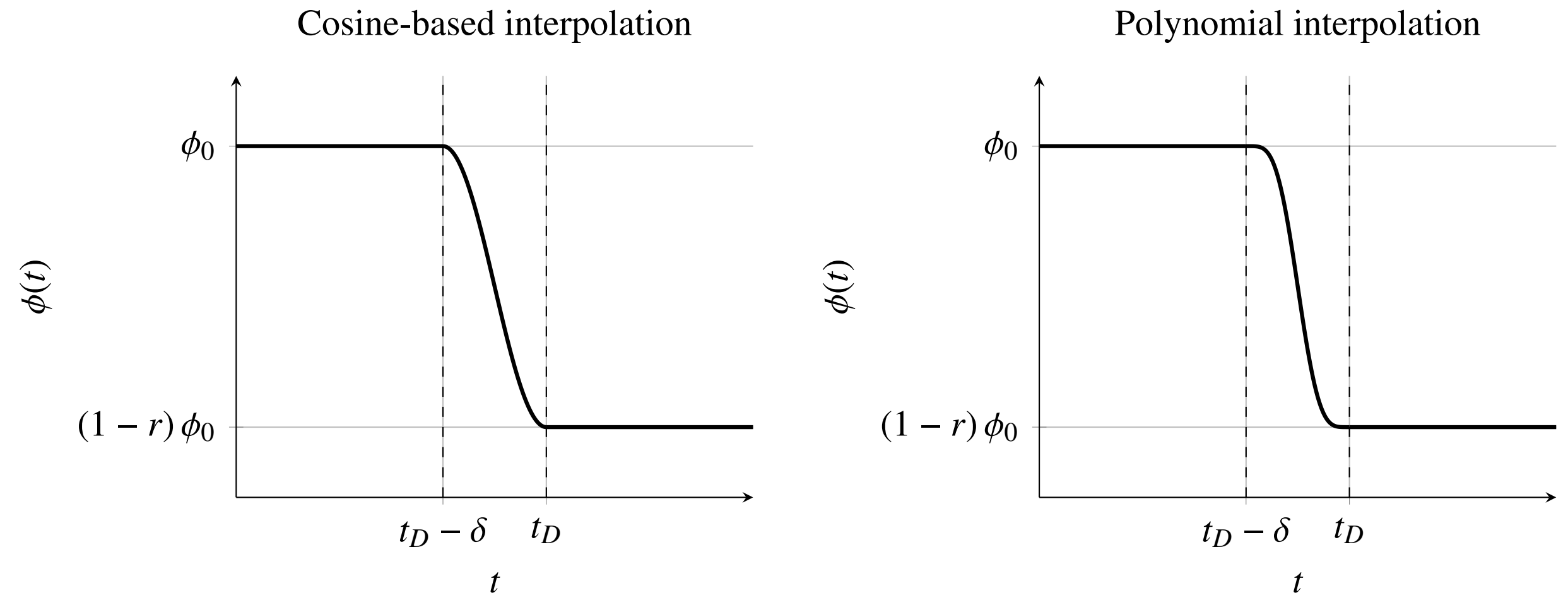

To handle such discrete policy changes accordingly close to reality, we simulate contact rates that transition over a small finite interval. More precisely, for a change effective at $t_D\in\mathbb{R}$ that changes a prior contact rate $\phi_0>0$ to $(1-r)\,\phi_0$, we define
\begin{align}\label{eq:phi}
	{\phi}(t) := 
	\begin{cases}
		\phi_{0},              & t \leq t_D  - \delta\\
		\widetilde{\phi}^*(t), & t \in (t_D - \delta, t_D) \\
		(1-r)\,\phi_{0} , & t \geq t_D. 
	\end{cases}
\end{align}
Here, $r\in(-\infty,1]$ denotes the relative change of the contact rate. Positive values of $r$ correspond to a contact reduction in which case we sometimes refer to the new contact rate as \textit{damped} value; negative values of $r$ allow for an increase beyond the baseline. The term $\delta>0$ defines the length of the transition interval, and $\widetilde{\phi}^*:(t_D - \delta,t_D)\rightarrow((1-r)\,\phi_{0}, \phi_{0})$ represents a transition contact rate to ensure a smooth transition; see, e.g.,~\cite{schmidt_gnn_2025}. Note that the length of the transition interval can be chosen arbitrarily for $0<\delta<\infty$. In~\cite{wendler_nonstandard_2026} and the MEmilio framework~\cite{bicker2026memilio} itself, we used the smoothed contact rates $\widetilde{\phi}^*=\widetilde{\phi}^C$ defined as
\begin{align}\label{eq:smooth_cos}
	\widetilde{\phi}^C(t):=\frac{1}{2}(r\phi_0)\cos\left(\pi\left(\frac{t-t_D+ \delta}{\delta}\right)\right)+\frac{1}{2}(2\phi_0-r\phi_0);
\end{align}
cf.~\cref{fig:smoother} (left). Based on the cosine function, it is guaranteed that ${\phi}\in\mathcal{C}^1(\mathbb{R})$. However, to fully benefit from the higher-order numerical schemes discussed in this paper, we need an interpolation function that satisfies ${\phi}\in\mathcal{C}^k(\mathbb{R})$ with $k>1$. Choosing the polynomial interpolation or \textit{smoothstep} function $\widetilde{\phi}^*=\widetilde{\phi}^{P_4}$ defined on $(t_D - \delta,t_D)$ as
\begin{align}\label{eq:smooth_pol}
	\begin{aligned}
		\widetilde{\phi}^{P_4}(t) &:= \phi_0 + ((1-r)\,\phi_0 - \phi_0) \left(126\left(\frac{t - t_D + \delta}{\delta}\right)^5 -420\left(\frac{t - t_D + \delta}{\delta}\right)^6 + 540 \left(\frac{t - t_D + \delta}{\delta}\right)^7\right.\\
		&\hspace*{4cm}\left.-315\left(\frac{t - t_D + \delta}{\delta}\right)^8 +70\left(\frac{t - t_D + \delta}{\delta}\right)^9\right)
	\end{aligned}
\end{align}
is sufficient to ensure ${\phi}\in\mathcal{C}^4(\mathbb{R})$; cf.~\cref{fig:smoother} (right) for a visualization.

\section{A discretized age-of-infection model for the solution of all compartments}\label{sec:discretization}

\subsection{Preliminaries}

We discretize the age-of-infection model derived in~\cref{sec:continuous_model} by applying a numerical scheme to the continuous model. Hereby, we distinguish between the solution scheme for the equation for the susceptibles, which mainly follows the ideas of~\cite{messina_positive_2022}, and the solution scheme for the numbers of infected and recovered individuals. In the latter, we suggest to use an additional finite difference (FD) scheme to approximate $S'(t)$ to reduce the number of function evaluations compared to approximating the integral within the analytical derivative, cf.~\eqref{eq:Sprime_explicit}.
For a clear distinction between continuous and discretized function definitions, we use the $\widehat{u}$ notation for the respective discretized version of a continuous function $u$. 

For a constant step size $\Delta t>0$, let us define a uniform discretization of $\mathbb{R}$ through $t_n:=n\Delta t$ for $n\in\mathbb{Z}$. Without loss of generality, we assume that there exists a nonpositive value $a\in \mathbb{Z}_{\leq 0}$ such that 
\begin{align*}
	t_{a} := a\,\Delta t = \widetilde{T}.
\end{align*}
As the (numerical) computation of~\eqref{eq:I_R_cont_finite} and~\eqref{eq:final_S_cont_finite} always involves past values on $[\widetilde{T},t_n]$, or on a numerically truncated subinterval if the corresponding kernels are negligible beyond a finite infection age, we cannot trivially adapt the step size on a considered interval. This fact renders the use of, e.g., Gaussian or adaptive Gau\ss-Kronrod quadrature less attractive, as these methods would require interpolation or additional evaluations of the solution history. As we similarly do not possess information on higher order derivatives, we cannot use Euler–Maclaurin's formula. As suggested by~\cite{messina_positive_2022}, we will use Gregory's quadrature rule. With Gregory's rule being based on the composite trapezoidal rule and adequately weighted end-point corrections, it satisfies our requirements. 

\subsection{Discretization of \texorpdfstring{$S$}{S}}
Analogously to~\cite{messina_positive_2022}, we discretize~\eqref{eq:final_S_cont_finite} by a direct quadrature method with quadrature weights of the $n_0$-th Gregory rule and obtain 
\begin{align}\label{eq:Sdiscret}
	\begin{aligned}
		\widehat{S}(t_n) = \widehat{S}(t_{a})\,\exp\Bigg(&\Delta t\,\Bigg[\sum_{j=a}^{n} \omega_{n-a,j-a}^{n_0} \,\rho(t_{n-j})\,\xi(t_{n-j})\,\gamma_I^R(t_{n-j})\left(\frac{\phi(t_n)}{N}\,\widehat{S}(t_j)+\frac{\phi(t_{n-j+a})}{N}\left(\widehat{R}(t_a)-N\right)\right) \Bigg] \\
		&- \Delta t^2\, \Bigg[\sum_{j=a}^n \omega_{n-a,j-a}^{n_0} \, \frac{\phi'(t_j)}{N}\sum_{k=a}^j \omega_{j-a, k-a}^{n_0}\rho(t_{j-k})\,\xi(t_{j-k})\,\gamma_I^R(t_{j-k})\, \widehat{S}(t_k)\Bigg]\Bigg)
	\end{aligned}
\end{align}
where the weights $\omega_{n-a,j-a}^{n_0}$ refer to Gregory rule $n_0$ and are defined as in~\cite{brunnerNumericalSolutionVolterra1986}; see also~\ref{app:gregory_weights} for a more detailed description of the Gregory weights. Here, the second sum appears due to our time-dependent contact rate which was not present in~\cite{messina_non-standard_2022}. Note that this defines an implicit equation which we solve by a fixed-point iteration. As an initial guess for $\widehat{S}(t_n)$, we use $\widehat{S}(t_{n-1})$ obtained at the previous time step. The iteration is terminated once the absolute difference between two successive approximations falls below $tol_{\text{FP}}=10^{-14}$.

As demonstrated by~\cite{messina_positive_2022} for the model with constant contact rate, we target convergence order $n_0+1$ when using quadrature weights of Gregory rule $n_0$ for the simplified formulation of the equation of the susceptible individuals. As a subtle requirement, we need at least $n_0$ prior values to apply Gregory rule $n_0$. 

\subsection{Discretization of \texorpdfstring{$I$}{I} and \texorpdfstring{$R$}{R}}
For the discretization of the equations $I(t)$ and $R(t)$, we require an approximation of the derivative ${S}'$. In~\eqref{eq:Sprime_explicit}, we derived an analytical expression for the derivative, making it possible to approximate the integral term with higher order quadrature rules. However, given a solution for $S(t)$, we can avoid evaluating this integral and substantially reduce the number of function evaluations or accesses to function values in the memory or cache. Using a FD scheme for approximating the derivative, only the five most recent evaluations are needed to compute a fourth order accurate approximation. For the sake of simplicity, we employ by default a fourth order FD scheme to ensure up to fourth order convergence of the numerical solution scheme. More precisely, we employ the backwards FD scheme 
\begin{align}\label{eq:FD4}
	g_4(t,f,\Delta t) := \frac{25\,{f}(t) -48\,{f}(t-\Delta t)+36\,{f}(t-2\Delta t)-16\,{f}(t-3\Delta t)+3\,{f}(t-4\Delta t)}{12\,\Delta t}.
\end{align}
and set
\begin{align}\label{eq:S_deriv_discrete}
	\widehat{S}'(t_j) =  g_4(t_j,\widehat{S},\Delta t)
\end{align}
for $j>a$. For sufficiently smooth functions, we can easily show that 
\begin{align}
	f'(t_j) - g_4(t_j,f,\Delta t) = \mathcal{O}(\Delta t^4).
\end{align}
For a general description of how to derive FD approximations of arbitrary order, we refer to~\cite[Chapter~1]{leveque_fd_schemes_2007}.

In additional numerical considerations, we will also briefly consider lower order FD schemes which are sufficient when utilizing lower order Gregory's rules. Here, we use
\begin{align}\label{eq:fd123}
	\begin{aligned}
		g_1(t,f,\Delta t)&:=\frac{{f}(t) -{f}(t-\Delta t)}{\Delta t},\\
		g_2(t,f,\Delta t)&:=\frac{3\,{f}(t) -4\,{f}(t-\Delta t)+\,{f}(t-2\Delta t)}{2\,\Delta t},\\
		g_3(t,f,\Delta t)&:=\frac{11\,{f}(t) -18\,{f}(t-\Delta t)+9\,{f}(t-2\Delta t)-2\,{f}(t-3\Delta t)}{6\,\Delta t},\\
	\end{aligned}
\end{align}
where
\begin{align}
	f'(t) - g_{m_0}(t,f,\Delta t) = \mathcal{O}(\Delta t^{m_0}),\quad m_0\in\{1,2,3\}.
\end{align}
With this, we can proceed to discretize the equations for $I$ and $R$, again using Gregory quadrature weights and obtain
\begin{align}\label{eq:I_R_discrete}
	\begin{aligned}
		\widehat{I}(t_n) &= \gamma_I^R(t_{n-a})\,\widehat{I}(t_{a}) -\Delta t \sum_{j=a}^{n} \omega_{n-a,j-a}^{n_0}\,\gamma_I^R(t_{n-j})\,\widehat{S}'(t_j) , \\
		\widehat{R}(t_n) &= N - \widehat{S}(t_n) - \widehat{I}(t_n).
	\end{aligned}
\end{align}  

\subsection{Composition of the error}\label{sec:error}

We distinguish the error introduced by time discretization from the error associated with the finite-history approximation. Let $u$ denote the solution of~\eqref{eq:S_cont_infty} and~\eqref{eq:I_R_cont_infty}, $u_{\widetilde{T}}$ the continuous finite-history solution of~\eqref{eq:final_S_cont_finite}~and~\eqref{eq:I_R_cont_finite}, and $\widehat{u}_{\Delta t,\widetilde{T}}$ the numerical solution of~\eqref{eq:Sdiscret}~and~\eqref{eq:I_R_discrete} with time step $\Delta t$. For each grid point $t_j$, $j\geq a$, the difference from the full-history solution decomposes as
\begin{align}\label{eq:error_decomposition_history}
	\widehat{u}_{\Delta t,\widetilde{T}}(t_j)-u(t_j)
	= \underbrace{u_{\widetilde{T}}(t_j)-u(t_j)}_{\text{history approximation error}}+\hspace{0.6cm}\underbrace{\widehat{u}_{\Delta t,\widetilde{T}}(t_j)-u_{\widetilde{T}}(t_j)}_{\text{numerical error}}.
\end{align}
For arbitrary but fixed $\widetilde{T}$, the first term is driven by the error that may be induced by~\eqref{eq:history_error_Sprime} and convergence of the numerical scheme concerns the second term. Consequently, comparison with a full-history reference solution measures both contributions. If the history approximation error is nonzero, refinement of the time step alone cannot produce convergence to that reference solution, although the numerical convergence order may be visible while the numerical error dominates. In the special case of exponential survival and age-independent $\rho\xi$, \cref{lemma:error} shows that the history related error vanishes.

In applications and, in particular, discretized settings, we have a $\widehat{T}$ which corresponds to the introduction of the first case(s). If the time of the first infection is known, $\widetilde{T}$ can be set to this value. Alternatively, the choice of $\widetilde{T}$ can be considered a modeling assumption that may introduce an approximation error if $\widetilde{T}>\widehat{T}$. We will examine the impact of different modeling errors in~\cref{sec:impact_epidemic}.

\subsection{Required function evaluations and initial step sizes}\label{sec:initial}

A particular problem of the numerical scheme is visible from~\eqref{eq:Sdiscret},~\eqref{eq:S_deriv_discrete}, and~\eqref{eq:I_R_discrete}: For Gregory rule $n_0$ with FD scheme of order $m_0$, we have insufficient past values for $n<a+\max\{n_0,m_0\}$ and the corresponding formulas can thus not be applied. 

As we require at least $n_0$ initial values for employing the $n_0$-th Gregory rule in~\eqref{eq:Sdiscret}, we increase the Gregory rule stepwise until we have approximated enough values to apply the $n_0$-th Gregory rule. This means that we lose accuracy over the first steps as the convergence order will not exceed $\mathcal{O}(\Delta t^{q+1})$ at $t_{a+q}$, $q=1,\ldots,n_0-1$ for $n_0>1$. For $n_0=1$, we only require one initial value, thus the stepwise increase of Gregory rules is not necessary. 
%
% To ensure that the largest error on the initial interval is $c_1\,\Delta t^{n_0+1}$, we need to adapt the step size for the initial steps of the model such that for an initial $\Delta t_{0}$, we have
% \begin{align}\label{eq:reducedh}
	% c_0\Delta t_{0}^2 \leq  c_1\Delta t^{n_0+1}\quad\text{and}\quad k\,\Delta t_{0} = \Delta t,\; k\in\mathbb{N}.
	% \end{align}
% Here, $c_0$ and $c_1$ are unknown constants and the first condition ensures the error to be driven by the $c_1\Delta t^{n_0+1}$. The second condition ensures that the progression with the first fine steps eventually reaches the coarse step $\Delta t$ exactly. Neglecting the constants, we can choose
% \begin{align}\label{eq:kh0}
	% k:=\lceil{\Delta t^{-\frac{n_0-1}{2}}}\rceil.
	% \end{align}

% {\color{blue}
	% SUGGESTION:
	
	To prevent the lower-order error introduced during initialization from asymptotically dominating the error of the higher-order scheme, we choose $\Delta t_0$ such that
	\begin{align}\label{eq:reducedh}
		c_0\Delta t_0^2 \leq c_1\Delta t^{n_0+1},
	\end{align}
	where $c_0$ and $c_1$ denote problem-dependent constants. This yields
	\begin{align}
		\Delta t_0
		\leq
		\sqrt{\frac{c_1}{c_0}}
		\Delta t^{(n_0+1)/2}.
	\end{align}
	Writing $\Delta t_0=\Delta t/k$, we obtain
	\begin{align}
		k\geq
		\sqrt{\frac{c_0}{c_1}}
		\Delta t^{-(n_0-1)/2}.
	\end{align}
	Since the constants $c_0$ and $c_1$ are generally unknown, we use the asymptotic scaling suggested by this relation and choose
	\begin{align}
		k=\left\lceil\Delta t^{-(n_0-1)/2}\right\rceil.
		\label{eq:kh0}
	\end{align}
	Thus, the choice of $\Delta t_0$ is intended to provide the appropriate asymptotic scaling of the initialization error rather than a strict bound on its magnitude.
	% }

After having computed $\widehat{S}$, we can compute $\widehat{S}'$. As the computations can be decoupled, we do not increase the order of the FD scheme step by step as done for the Gregory rule. Instead, we use a forward FD scheme of $m_0$-th order to compute the derivative for $t_j\in\{t_{a},\dots,t_{a+m_0-1}\}$ before switching to a backwards FD scheme for $t_j>t_{a+m_0-1}$. Note that, here, we have to ensure to compute enough values for $S(t)$ beforehand.

For the forward FD approximation we use
\begin{align}\label{eq:FD_forw}
	\begin{aligned}
		h_1(t,f,\Delta t)&:=\frac{-{f}(t) +{f}(t+\Delta t)}{\Delta t},\\
		h_2(t,f,\Delta t)&:=\frac{-3\,{f}(t) +4\,{f}(t+\Delta t)-\,{f}(t+2\Delta t)}{2\,\Delta t},\\
		h_3(t,f,\Delta t)&:=\frac{-11\,{f}(t) +18\,{f}(t+\Delta t)-9\,{f}(t+2\Delta t)+2\,{f}(t+3\Delta t)}{6\,\Delta t},\\
		h_4(t,f,\Delta t) &:=\frac{-25\,{f}(t) +48\,{f}(t+\Delta t)-36\,{f}(t+2\Delta t)+16\,{f}(t+3\Delta t)-3\,{f}(t+4\Delta t)}{12\,\Delta t},
	\end{aligned}
\end{align}
where
\begin{align}
	f'(t) - h_{m_0}(t,f,\Delta t) = \mathcal{O}(\Delta t^{m_0}),\quad m_0\in\{1,2,3,4\}.
\end{align}
With this, we set
\begin{align}\label{eq:Sderiv_FD_forw}
	\widehat{S}'(t_j) =  h_{m_0}(t_j,\widehat{S},\Delta t)
\end{align}
for $t_j\in\{t_{a},\dots,t_{a+{m_0}-1}\}$.

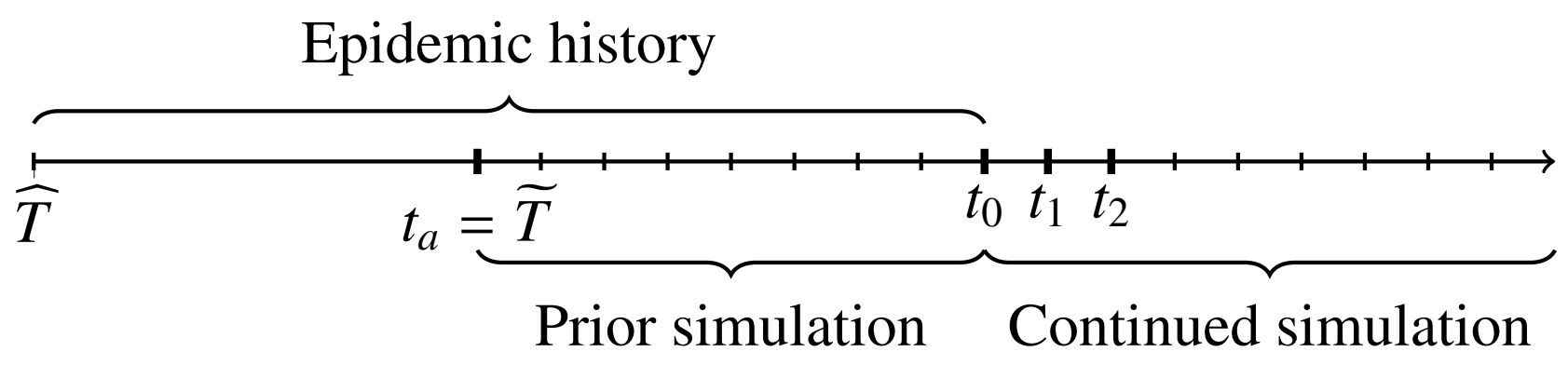
\begin{figure}
	\centering
	\begin{tikzpicture}[scale=0.45,every node/.style={
			scale=0.85, font=\Large}]
		
		% Axis
		\draw[thick,->] (12,0) -- (36,0);
		
		% equidistant ticks
		\draw[thick] (12,-0.14) -- (12,0.14);
		\foreach \x in {19,...,35}{
			\draw[thick] (\x,-0.14) -- (\x,0.14);
		}
		\draw[ultra thick] (12,-0.2)  -- (12,-0.2);
		\draw[ultra thick] (19,-0.2) -- (19,0.2);
		\draw[ultra thick] (27,-0.2) -- (27,0.2);
		\draw[ultra thick] (28,-0.2) -- (28,0.2);
		\draw[ultra thick] (29,-0.2) -- (29,0.2);
		
		% Labels: tick number = coordinate + 1
		\node[below] at (12,-0.1)  {$\widehat{T}$};
		\node[below] at (19,-0.1) {$t_a=\widetilde{T}$};
		\node[below] at (27,-0.1) {$t_0$};
		\node[below] at (28,-0.1) {$t_1$};
		\node[below] at (29,-0.1) {$t_2$};
		
		% Brace above ticks 1--28
		\draw[thick,decorate,decoration={brace,amplitude=5pt}]
		(12,0.6) -- (27,0.6);
		\node[above] at (19.5,1.1) {Epidemic history};
		
		\draw[thick,decorate,decoration={brace,amplitude=5pt}]
		(27,-1.4) -- (19,-1.4);
		\node[below] at (23,-2) {Prior simulation};     
		
		\draw[thick,decorate,decoration={brace,amplitude=5pt}]
		(36,-1.4) -- (27,-1.4);
		\node[below] at (31.5,-2) {Continued simulation};  
		
	\end{tikzpicture}
	\caption{\textbf{Epidemic history, model consideration, and simulation windows.}}
	\label{fig:restart}
\end{figure}

These considerations are highly relevant for the first time steps to accurately approximate the right solution. However, in practical situations where our model might be iteratively advanced and we might already have historical data, i.e., the simulation with the model actually corresponds to a restart of a prior simulation until $t_0$; cf.~\cref{fig:restart}. No evaluations points are available for $t<\widetilde{T}$ but at $t=t_0$ we actually restart the model simulation with historical data.

\subsection{Implementation} 

The prior considerations yield the straightforward implementation of the numerical scheme as depicted in~\cref{alg:ide_solver}. The set of input parameters are defined via
\begin{align}\label{eq:alg_params}
	\mathcal{P} := \{\gamma_I^R, \rho, \xi, \phi,\phi', n_0, m_0, tol_{\text{FP}},\widehat{S}(t_{a}), \widehat{I}(t_{a}), \widehat{R}(t_{a})\}.
\end{align}
To reduce the error in the initial interval, we can apply our numerical scheme with a smaller initial step size $\Delta t_0$ up to $t_{switch}$, cf.~\eqref{eq:reducedh}, before switching to the larger step size $\Delta t$ for the remaining simulation steps.

We start with the initial step size $\Delta t_0$. First, we compute all values $\widehat{S}(t_j), t_j\in\{t_{a},\dots,t_{switch}\}$. Then, we use these values to compute the derivative values $\widehat{S}'(t_j)$ together with the compartment values $\widehat{I}(t_j)$ and $\widehat{R}(t_j)$ until we reach $t_{switch}$. From there on, we use the step size $\Delta t$ and repeat the same steps for  $t_j\in\{t_{switch}+\Delta t,\dots,t_{\max}\}$
This is depicted in~\cref{alg:ide_solver}.
The implementation shown here is straightforward, in order to improve the efficiency and performance of~\cref{alg:ide_solver}, the main loops could be combined and rearranged.

\begin{algorithm}[!ht]
	\caption{Discrete scheme.}
	\label{alg:ide_solver}
	\begin{algorithmic}[1]
		\Require Step sizes $\Delta t$, time points $t_{a}$, $t_{\max}$, parameter set $\mathcal{P}$ according to~\eqref{eq:alg_params}
		\Statex \hspace{\algorithmicindent}where $p=\max\{n_0-1,2m_0-1\}$.
		\State $t_{\mathrm{switch}} \gets t_a + p\Delta t$
		\State $\Delta t_0 \gets \Delta t/k$ with~\eqref{eq:kh0}
		%\Statex
		\State $t_j \gets t_{a}+\Delta t_0$
		\While{$t_j \leq t_{\mathrm{switch}}$}
		\State Compute $\widehat{S}(t_j)$ according to~\eqref{eq:Sdiscret} with (increasing) Gregory rule $n_0$ and step size $\Delta t_0$.
		\State $t_j \gets t_j + \Delta t_0$
		\EndWhile
		\State Compute $\widehat{S}'(t_a)$ according to~\eqref{eq:Sderiv_FD_forw} and step size $\Delta t_0$.
		\State $t_j \gets t_{a}+\Delta t_0$
		\While{$t_j \leq t_{\mathrm{switch}}$}
		\State Compute $\widehat{S}'(t_j)$ according to~\eqref{eq:Sderiv_FD_forw} or~\eqref{eq:S_deriv_discrete} and step size $\Delta t_0$.
		\State Compute $\widehat{I}(t_j)$ and $\widehat{R}(t_j)$ according to~\eqref{eq:I_R_discrete} with (increasing) Gregory rule $n_0$ and step size $\Delta t_0$.
		\State $t_j \gets t_j + \Delta t_0$
		\EndWhile
		\Statex 
		\State $t_j \gets t_{\mathrm{switch}}+\Delta t$
		\While{$t_j \leq t_{\max}$}
		\State Compute $\widehat{S}(t_j)$ according to~\eqref{eq:Sdiscret} with Gregory rule $n_0$ and step size $\Delta t$.
		\State $t_j \gets t_j + \Delta t$
		\EndWhile
		
		\State $t_j \gets t_{\mathrm{switch}}+\Delta t$
		\While{$t_j \leq t_{\max}$}
		\State Compute $\widehat{S}'(t_j)$ according to~\eqref{eq:S_deriv_discrete} and step size $\Delta t$.
		\State Compute $\widehat{I}(t_j)$ and $\widehat{R}(t_j)$ according to~\eqref{eq:I_R_discrete} with Gregory rule $n_0$ and step size $\Delta t$.
		\State $t_j \gets t_j + \Delta t$
		\EndWhile
	\end{algorithmic}
\end{algorithm}

\subsection{Observations regarding the numerical solution} 
We observe that the numerical scheme as in~\cref{sec:discretization} is mass-conserving which follows directly from the discrete formulas for the compartments as in~\eqref{eq:Sdiscret} and~\eqref{eq:I_R_discrete}.
\begin{lemma}
	For all $j\in\mathbb{Z}$, $j\geq a$, it holds that $\widehat{S}(t_{j})+\widehat{I}(t_{j})+\widehat{R}(t_{j})=N$.
\end{lemma}
Furthermore, it is easy to show that the positivity of $\widehat{S}$ is conserved due to the reformulation with an exponential operator. 
\begin{lemma}
	Let $\widehat{S}(t_a)>0$. Then it holds that $\widehat{S}(t_{j})>0$ for all $j\in\mathbb{Z}$, $j>a$.
\end{lemma}

\section{Numerical results}\label{sec:numerical_results}
The presented IDE model together with our numerical experiments were implemented in C++ as a part of our high performance modular epidemics
simulation software MEmilio~\cite{bicker2026memilio}. We present both qualitative and detailed quantitative results on the model and its numerical solution through higher order schemes. We consider different settings regarding Gregory rules and FD approximations as well as different smooth interpolation functions on contact change points. Moreover, we examine the impact of the initialization horizon on the epidemic outcomes.

\subsection{Analysis of the convergence order}

In this section, we provide the numerical convergence order for our scheme as introduced in~\cref{sec:discretization}. Starting with a constant contact rate, we first consider the convergence order for different Gregory rules combined with a FD approximation of order 4. Second, we provide the convergence orders when using different FD approximation orders. Third, we consider a time-dependent contact rate. 

To approximate the convergence order, we need to compute a ground truth solution. Instead of solving the IDE formulation directly, we consider an exponential sojourn time distribution. We make use of~\cref{lemma:ode,lemma:error} to apply a high order Runge--Kutta scheme that delivers the ground truth as the solution of the corresponding ODE-based initial value problem (IVP). 
We denote the corresponding solution by $u_{\text{ODE}}^*$. The results of the IDE model are computed with different time steps and then compared to the ground truth. 

The relative error is computed by normalizing the absolute error with the ground truth norm, i.e.,
\begin{align}\label{eq:max1}
	err_{\text{rel}} := \frac{\left\lVert\widehat{u}_{\text{IDE}}- u_{\text{ODE}}^*\right\rVert_\infty}{\left\lVert u_{\text{ODE}}^*\right\rVert_\infty}
\end{align}
for any compartment $\widehat{u}_{\text{IDE}}(t)$ representing one element in $\left\{\widehat{S}(t),\, \widehat{I}(t),\, \widehat{R}(t)\right\}$. The maximum norm is defined on the simulation interval through
\begin{align}\label{eq:max2}
	\left\lVert u \right\rVert_\infty := \max_{j=0,\dots,n} \left\vert u(t_j) \right\vert.
\end{align}

We use an exponential sojourn time distribution with mean $T_I=2$. We consider an initial point in time $\widehat{T}=-50$ and start with an initial population of 
\begin{align*}
	\widehat{S}(\widehat{T})=9\,999\,000, \quad \widehat{I}(\widehat{T})=1\,000,\quad \widehat{R}(\widehat{T})=0.
\end{align*}
We consider the ODE-based IVP with $t_{0,\text{ODE}}=\widehat{T}$ and solve until $t_\text{max}=100$. The IVP is solved with high precision as obtained through a Runge-Kutta solver of eighth order with constant step size $\Delta t = 10^{-6}$. The corresponding outcome is $u_{\text{ODE}}^*$.

In order to avoid the issue around the initial time point as described in~\cref{sec:discretization} and~\cref{alg:ide_solver}, we consider a restarted simulation as depicted in~\cref{fig:restart}. We set $\widetilde{T}=-40$ and thereby make implicit use of~\cref{lemma:ode,lemma:error}, i.e., that in case of exponential survival, the IDE model converges against the ODE model, independent of the choice of $\widetilde{T}$. In particular, we use the values for the compartments obtained via the ODE model, set them as high order approximations to the true values of the IDE model, and also use these to compute the derivative $S'$ by applying a FD scheme until $t_0=0$. Then, we simulate the IDE model from $t_{0,\text{IDE}}=0$ to $t_\text{max}$ and compute the error only over the simulation horizon. 

Further model parameters are set to $\rho(t)=0.8$, $\xi(t)=1$ and $\phi(t)=0.73$.

\subsubsection{Varying Gregory rules for constant contact rate}\label{sec:baseline}

First, we fix the FD approximations for $\widehat{S}'(t)$ to 4th order and consider the convergence of the numerical schemes with different Gregory rules of order $n_0\in\{1,2,3\}$ and $\Delta t\in\{10^0,10^{-1},10^{-2},10^{-3}\}$.

\begin{figure}[!ht]
	\centering
	\includegraphics[width=0.9\linewidth]{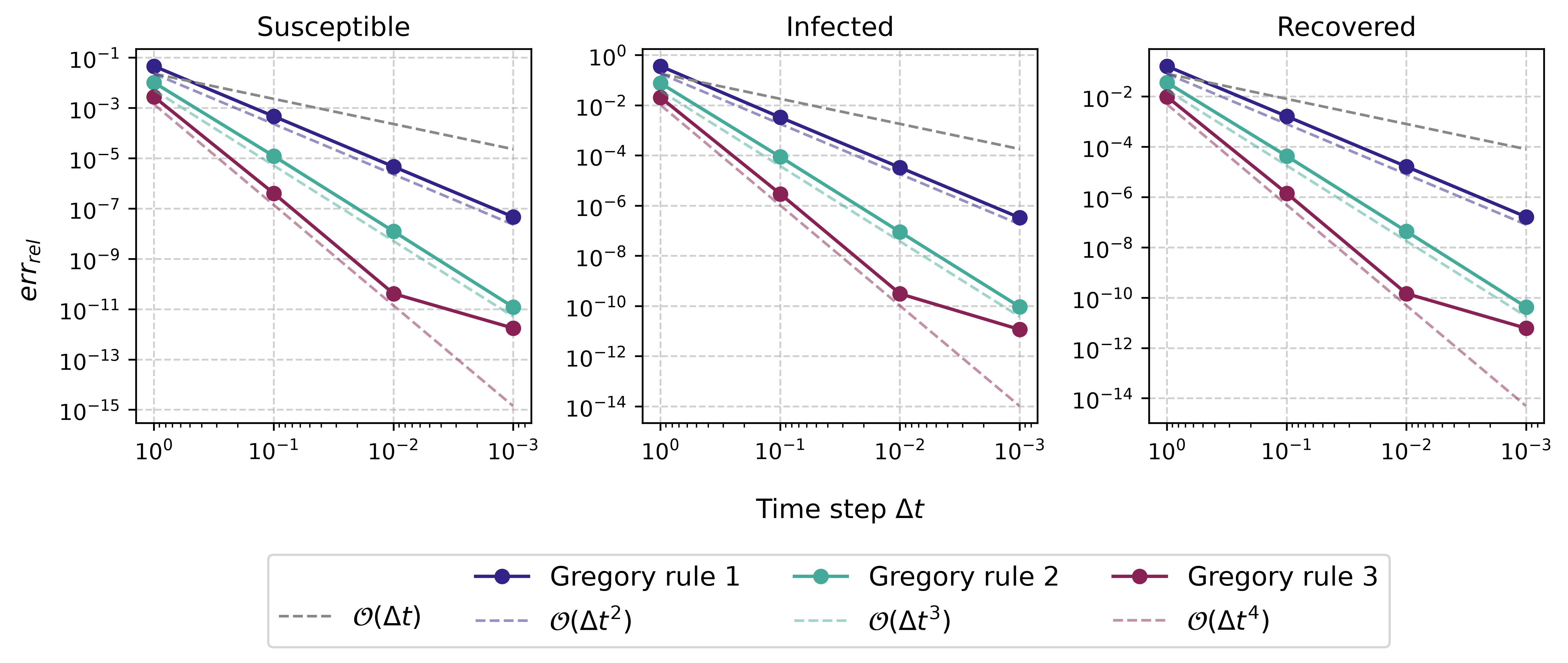}
	\caption{\textbf{Convergence of the suggested numerical schemes with different Gregory rules and 4th order FD approximation.} Convergence is computed and shown for the number of susceptibles (left), infected (center) and recovered (right) individuals using the maximum norm of the relative error as given through~\eqref{eq:max1},~\eqref{eq:max2}.}
	\label{fig:baseline}
	% Results from 2026-10-11
\end{figure}

From~\cref{fig:baseline}, we observe orders of convergence of $n_0+1$ for Gregory rules $n_0\in\{1,2\}$ in all compartments. The same holds for $n_0=3$ and $\Delta t\in\{10^0,10^{-1},10^{-2}\}$. For the smallest considered time step of $\Delta t=10^{-3}$, the error of Gregory rule 3 deviates from the optimal convergence behavior in all compartments. However, we note that the extrapolated relative error is close to machine precision for all compartments (i.e., $10^{-15}$ to $10^{-14}$) and hence we summarize that the baseline scheme is capable of attaining numerical convergence orders as intended.

\subsubsection{Impact of varying FD approximations}

In~\eqref{eq:FD4}, we have deliberately chosen more than the required FD approximation order to obtain the intended convergence order for different Gregory rules. In the following, we show how the convergence order is affected by the FD scheme. To do so, we consider FD approximations of order one, two and three as given in~\eqref{eq:fd123} to approximate $\widehat{S}'$. All other parameters follow the previous setting.

\begin{figure}[!ht]
	\centering
	\includegraphics[width=0.9\linewidth]{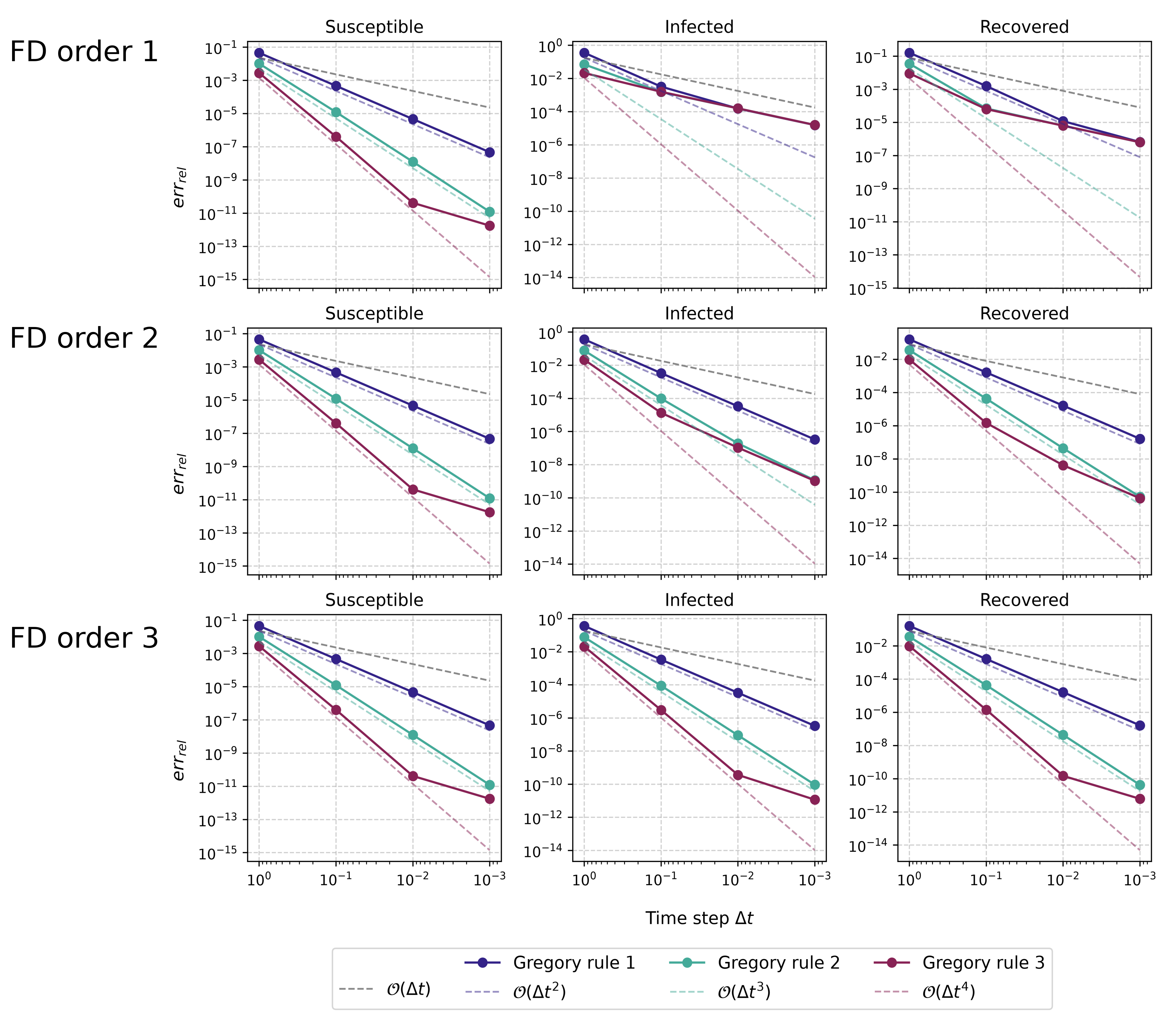}
	\caption{\textbf{Convergence of the numerical schemes with different Gregory rules for different FD approximations.} Different rows depict different FD orders. In each row, convergence is computed and shown for the number of susceptibles (left), infected (center) and recovered (right) individuals using the maximum norm of the relative error as given through~\eqref{eq:max1},~\eqref{eq:max2}.}
	\label{fig:fdorders}
	% Results from 2026-10-11
\end{figure}

As expected, the choice of the FD scheme does not influence the convergence in the approximated number of susceptibles (cf. left panels in~\cref{fig:fdorders}); in contrast to the number of infected and recovered individuals (center and right panels), where the FD scheme is used. For FD orders 1 and 2, the convergence order is dominated by the limited minimum order of convergence of the Gregory quadrature scheme and the FD scheme. We observe linear convergence when using a FD order 1 and quadratic convergence for FD order 2 for all considered Gregory rules, see~\cref{fig:fdorders}. Accordingly, we would expect second-order convergence for Gregory rule 1 and third-order convergence for Gregory rules 2 and 3 when using 3rd order FD. However, as shown in~\cref{fig:fdorders}, we obtain almost identical convergence for Gregory rule 3 and FD order 3 as for FD order 4. We did not further investigate if this effect is obtained by an additional cancellation of highest order error.

%Overall, we state that the convergence order of the approximations in $I(t)$ and $R(t)$ depends both on the Gregory rule $n_0$ as well as on the order of the FD scheme $m_0$. Numerically, we observe
%\begin{align*}
%p\geq\min\{n_0+1,m_0\}\quad\text{and}\quad p\leq \min\{n_0+1,m_0+1\}
%\end{align*}
%for $I$ and $R$ unless floating-point errors dominate. 
%With this, we implicitly confirmed that~\cref{lemma:ode,lemma:error} hold.

Overall, the numerical results indicate that the convergence order of the approximations of $I(t)$ and $R(t)$ depends on both the Gregory rule $n_0$ and the order $m_0$ of the FD scheme. Across the considered test cases, and unless floating-point errors dominate, we observe
\begin{align*}
	p\geq\min\{n_0+1,m_0\}\quad\text{and}\quad
	p\leq\min\{n_0+1,m_0+1\}.
\end{align*}
In particular, some configurations exhibit an empirical convergence order of $m_0+1$, exceeding the nominal order $m_0$ of the FD approximation. We do not provide a theoretical explanation for this apparent superconvergence here. The numerical results are furthermore consistent with the analytical results in~\cref{lemma:ode,lemma:error}.

\subsubsection{Convergence orders for a time-dependent contact rate}

In the previous examinations, we restricted ourselves to a constant contact rate $\phi(t)=\phi$ and evaluated the individual impacts of Gregory rules and FD approximation to the convergence rate. We now consider time-dependent contact rates as allowed by our generalized model formulation~\eqref{eq:final_S_cont_finite}. We also compare different interpolation functions to obtain a continuous transition for jumps in the contact rates as introduced in~\cref{sec:continuous_model}. In particular, we consider two smoothing interpolations, the $\mathcal{C}^{1}$ cosine interpolation of~\eqref{eq:smooth_cos} as used by default in MEmilio and the $\mathcal{C}^{4}$ polynomial smoothstep function of~\eqref{eq:smooth_pol}.

In both scenarios, we start as in the previous scenarios with a constant contact rate of $\phi=0.73$ before applying a 10\% reduction, corresponding to $r=0.1$ and damping factor $1-r=0.9$, at time $t_D=3.372907$ with a transition interval of length two, i.e., we set $r=0.1$ and $\delta=2$, cf. the notation in~\eqref{eq:phi}. To avoid potential effects of the damping applied on the discretized grid, we chose $t_D$ such that for all considered step sizes we have $t_D\neq t_j, j\in\mathbb{Z}$. Here, we use the analytical derivative of the contact rate $\phi'$ but proceed with the FD approximation for computing $\widehat{S}'$ to investigate the influence of the interpolation functions on the previously considered numerical scheme. We reduce the simulation interval to $t_\text{max}=5$ while all other parameters follow the previous setting.

\paragraph{Cosine interpolation function}

The cosine interpolation function has been made available within the MEmilio framework~\cite{bicker2026memilio} our IDE-based model is implemented with.

\begin{figure}[!ht]
	\centering
	\begin{subfigure}[a]{0.9\linewidth}
		\centering
		\includegraphics[width=0.9\linewidth]{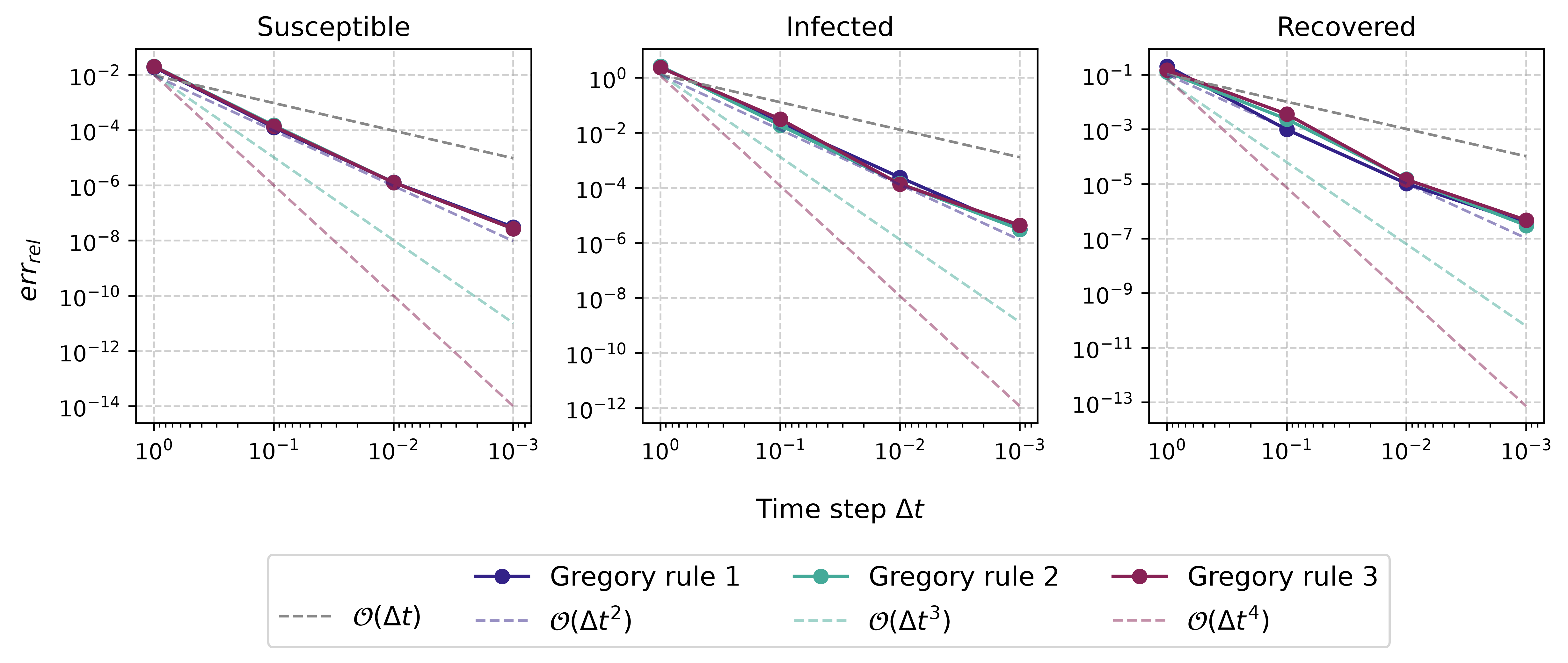}
	\end{subfigure}
	
	\begin{subfigure}[b]{0.9\linewidth}
		\centering
		\includegraphics[width=0.9\linewidth]{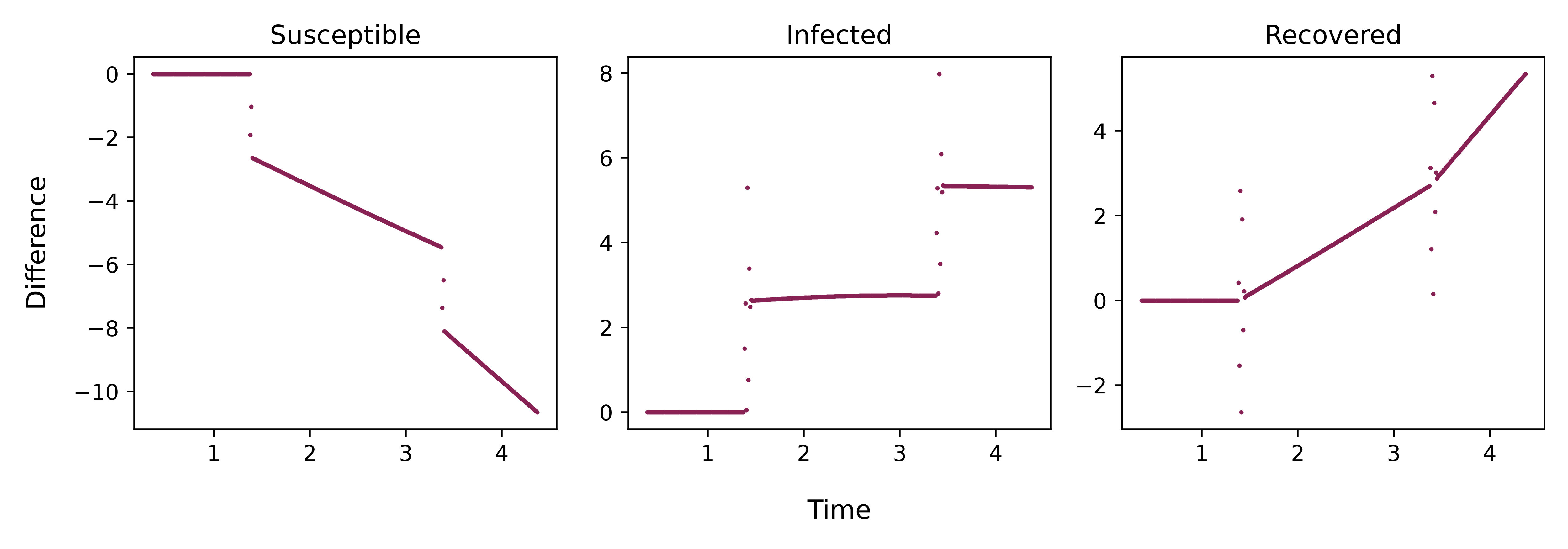}
	\end{subfigure}
	% Results from 2026-10-09
	\caption{\textbf{Convergence and difference per time step of the suggested numerical schemes with different Gregory rules, 4th order FD approximation, and cosine interpolation function for the time dependent contact rate.} Convergence is computed and shown for the number of susceptibles (top left), infected (top center) and recovered (top right) individuals using the maximum norm of the relative error as given through~\eqref{eq:max1},~\eqref{eq:max2}. Difference per time point between ground truth obtained with ODE model and simulation results of numerical scheme with 3rd Gregory rule, 4th order FD approximation and time step $\Delta t=10^{-2}$ for the number of susceptibles (bottom left), infected (bottom center) and recovered (bottom right).}
	\label{fig:smoothercosine}
\end{figure}

For all compartments and Gregory rules, we observe quadratic convergence, cf.~\cref{fig:smoothercosine} (top).
Due to reduced regularity of the integrands in~\eqref{eq:S_cont_infty_A} induced by the cosine interpolation function, we do not observe higher orders of convergence, even for Gregory rules 2 and 3. When considering the difference between the results of the ground truth and the simulation results of the IDE model for the susceptible individuals as depicted in the bottom row of~\cref{fig:smoothercosine}, we observe jumps around the beginning and end of the transition interval (at $t_j\in\{1.38, 1.39, 3.38, 3.39\}$), i.e., the first two discrete points after the start and end of the transition interval, respectively. When applying Gregory rule 3, we apply end point corrections to the first and last three summands which are computed via higher order finite difference operators; cf.~\cite[p.74]{brunnerNumericalSolutionVolterra1986}. In the case of $t_j\in\{1.38, 1.39, 3.38, 3.39\}$, these end point corrections are applied across start and end point of the transition interval where we have the reduced regularity of the cosine interpolation function. This leads to reduced orders of convergence as we assume higher regularity when deriving the correction terms. For the infected and recovered individuals, the difference between the ground truth and simulation results show even more jumps around the start and end of the transition interval. The jumps in the susceptibles are transferred via the FD approximation of the derivative $\widehat{S}'$ to the infected and recovered compartments.

\paragraph{Polynomial smoothstep function}

As the limited convergence order for the cosine interpolation function might come without surprise, we consider the polynomial smoothstep function which yields $\phi\in\mathcal{C}^4(\mathbb{R})$. 

\begin{figure}[!ht]
	\centering
	\begin{subfigure}[a]{0.9\linewidth}
		\centering
		\includegraphics[width=0.9\linewidth]{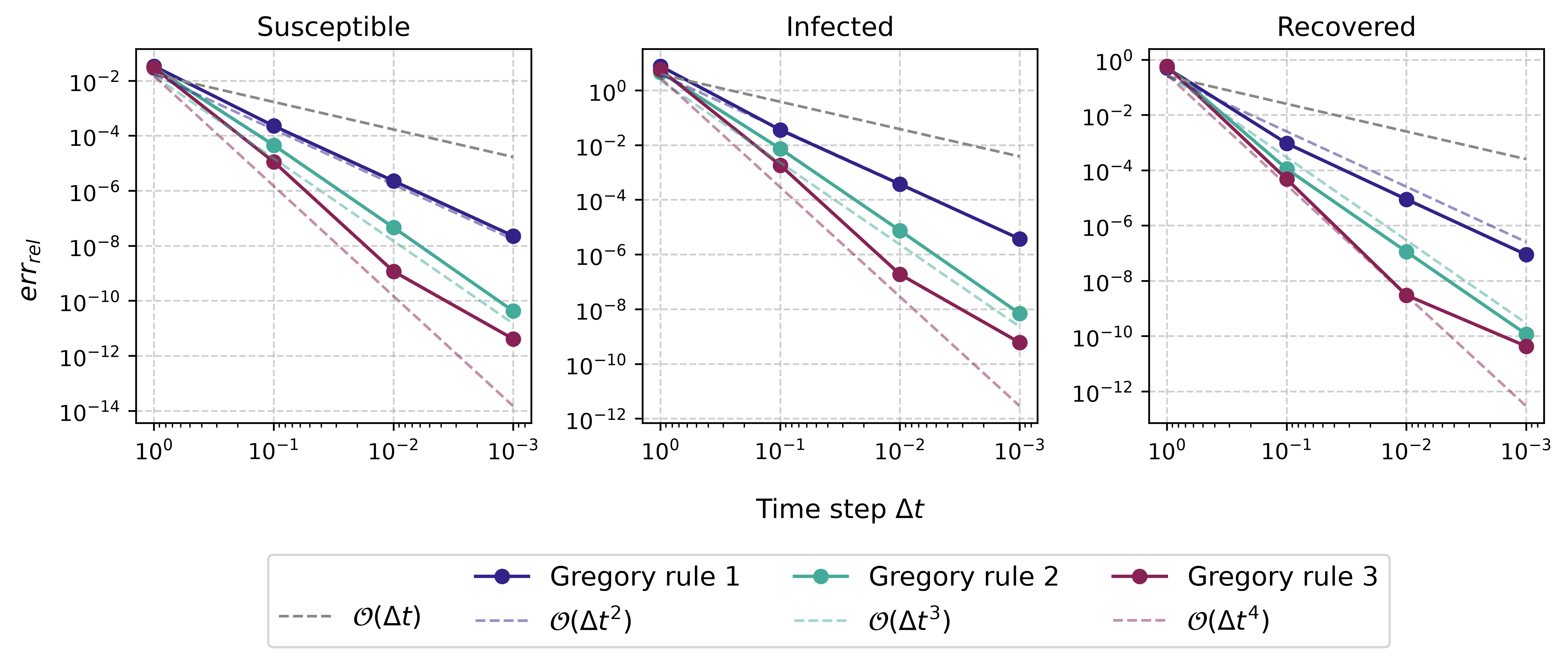}
	\end{subfigure}
	
	\begin{subfigure}[b]{0.9\linewidth}
		\centering
		\includegraphics[width=0.9\linewidth]{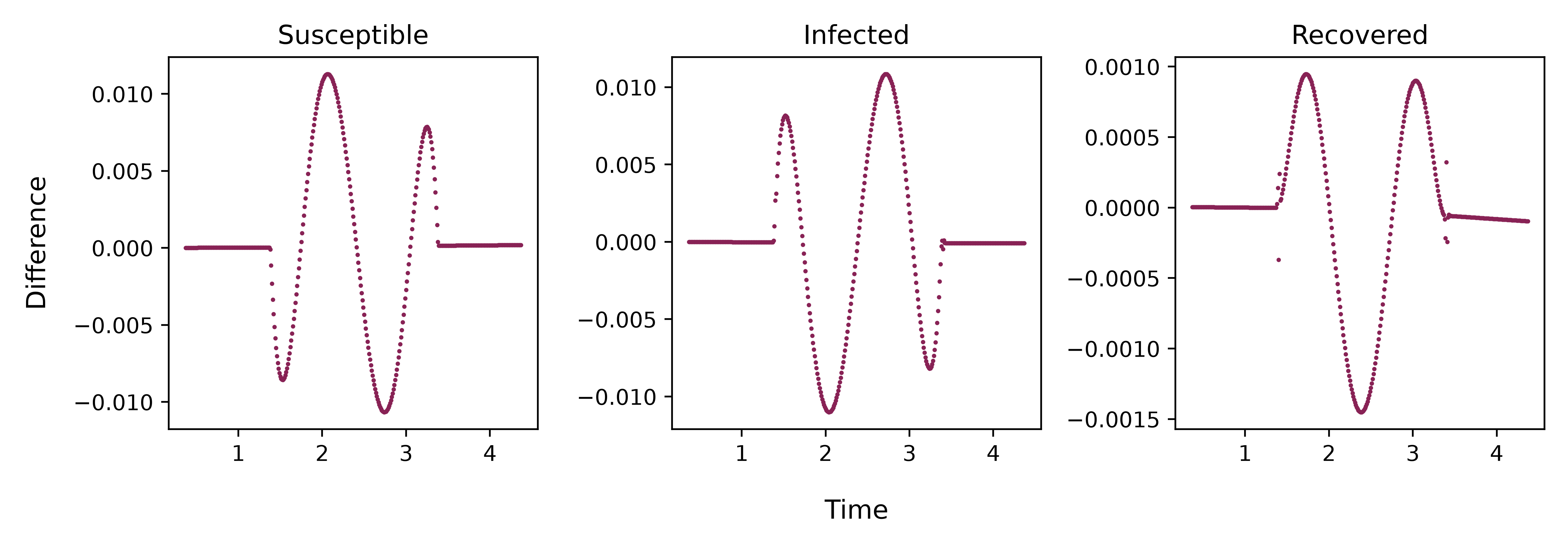}
	\end{subfigure}
	% Results from 2026-10-09
	\caption{\textbf{Convergence and difference per time step of the suggested numerical schemes with different Gregory rules, 4th order FD approximation, and $\mathcal{C}^4$ smoothstep function for the time dependent contact rate.} Convergence is computed and shown for the number of susceptibles (top left), infected (top center) and recovered (top right) individuals using the maximum norm of the relative error as given through~\eqref{eq:max1},~\eqref{eq:max2}. Difference per time point between ground truth obtained with ODE model and simulation results of numerical scheme with Gregory rule 3, 4th order FD approximation and time step $\Delta t=10^{-2}$  for the number of susceptibles (bottom left), infected (bottom center) and recovered (bottom right).}
	\label{fig:smoothstep}
\end{figure}

For all compartments, we observe convergence of $n_0+1$-th order for Gregory rule $n_0\in\{1,2,3\}$; cf.~\cref{fig:smoothstep} (top). When using the $\mathcal{C}^4$ smoothstep function for smoothing the jumps on some interval, we observe that the error per time step is the largest around the transition interval, fading out towards the interval ends; cf.~\cref{fig:smoothstep} (bottom). Since the contact rate $\phi$ is sufficiently smooth, we do not observe deteriorated errors around the end points of the transition interval for the susceptible individuals and, overall, obtain higher order convergence. For the infected and recovered individuals, we observe small jumps at around the start and end of the transition interval which are due to the kink in the difference of the susceptible individuals which is transferred to the subsequent compartments via its derivative $\widehat{S}'$. However, these errors do not limit the convergence.  

\paragraph{Smoothstep summary}

\begin{figure}[!ht]
	\centering
	\includegraphics[width=0.5\linewidth]{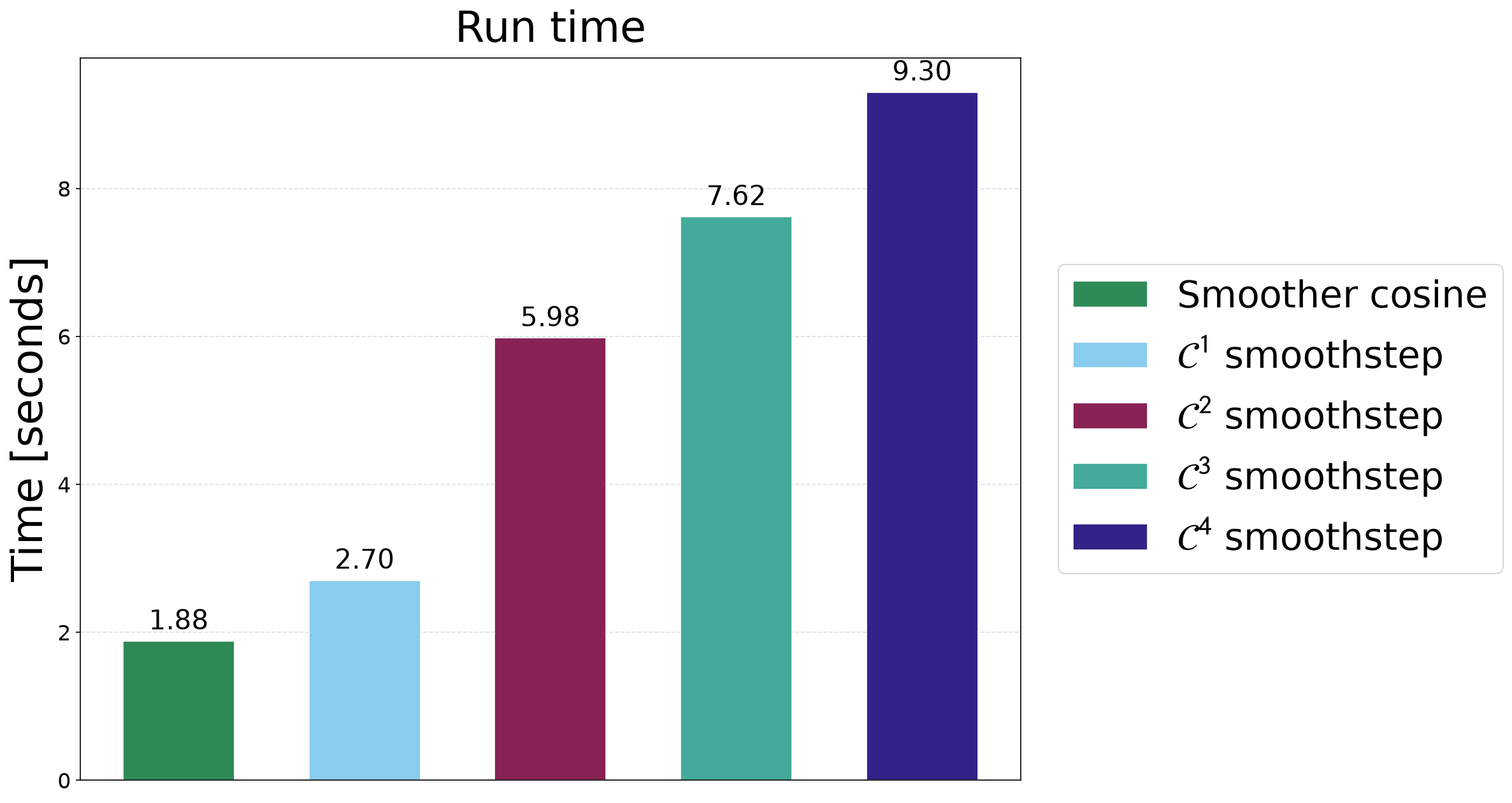}
	\caption{\textbf{Runtime comparison of the evaluation of cosine interpolation and different orders of polynomial smoothstep functions.} All functions smooth a damping from 1 to 0 on the transition interval $[0,100]$, see~\ref{app:smoothstep_def} for the definition of the polynomial smoothstep functions in $\mathcal{C}^k, k\in\{1,2,3\}$. The results were obtained by evaluating the corresponding functions on $1\,000\,001$ equidistant points on the considered interval. The figure shows the mean run time of 100 runs. Measurements were conducted on an Intel Xeon ‘‘Skylake’’ Gold 6132, 2.60 GHz with 14 cores per socket.}
	\label{fig:runtime_smoothstep}
\end{figure}

The formulation with smooth interpolation functions enables a practical setting or calibration of contact rates through the computation of a continuous solution on the entire simulation interval. In cases where higher order numerical schemes $n_0>2$ are required by the setting or the user, a $\mathcal{C}^k$ polynomial smoothstep function, with $k>2$, should be preferred. For cases where, however, quadratic convergence order is sufficient, the cosine interpolation function is, first, sufficient and, second, also more efficient than polynomial smoothstep functions; cf.~\cref{fig:runtime_smoothstep}.

\subsection{Impact of the included history on the convergence order}\label{sec:numerical_impact_init_conditions}

After having considered different numerical factors such as Gregory rule, FD approximation and interpolation of contact change points, we will now consider the convergence at the start of the simulation $\widetilde{T}$. For this, we will restrict ourselves again to a constant contact rate. As described in~\cref{sec:initial}, we require $n_0$ starting values to obtain convergence of order $n_0+1$ when applying Gregory rule $n_0$. In the previous sections, we have assumed that we already conducted a simulation with high accuracy for time points $t_j\in\{t_a,\dots,t_0\}$. This provided sufficiently accurate historic values to obtain the expected convergence order from the restart on. Now, we will investigate the case where we are not provided with suitable, previous simulation results and use a smaller step size at the start of the simulation to reduce the initial error below the higher numerical error at later time steps; cf. the discussion around~\eqref{eq:kh0}. 

\begin{figure}[!ht]
	\centering
	\includegraphics[width=0.9\linewidth]{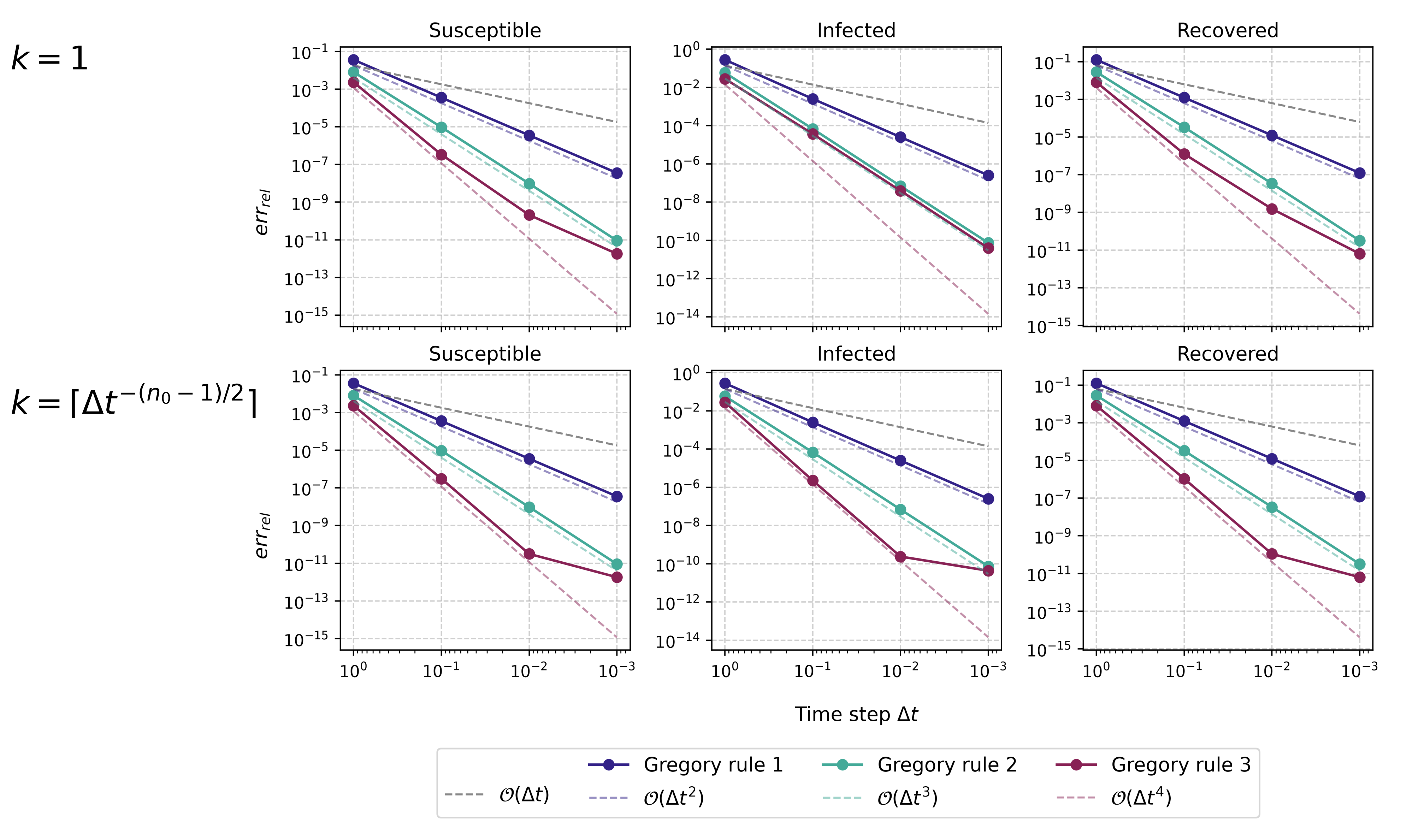}
	\caption{\textbf{Convergence of the numerical schemes with different Gregory rules without and with time-step adaptation.} The first row depicts the scenario without choosing a smaller initial step size such that $k=1$ while the bottom row shows the convergence when adapting the initial step size according to $\Delta t_0=\Delta t/k$ with $k=\left\lceil\Delta t^{-(n_0-1)/2}\right\rceil$, cf.~\eqref{eq:kh0}, until $t_{switch}$ before switching to the target step size $\Delta t$. In each row, convergence is computed and shown for the number of susceptibles (left), infected (center) and recovered (right) individuals using the maximum norm of the relative error as given through~\eqref{eq:max1},~\eqref{eq:max2}.}
	\label{fig:adapted_timestep}
	% Results from 2026-11-04
\end{figure}

We first show that the use of a constant step size $\Delta t$ from the start and over the whole horizon from $\widetilde{T}$ until $t_{\max}$ does indeed yield a reduced convergence order. In the second scenario, we use a finer step size $\Delta t_0$ on from $\widetilde{T}$ to $t_{switch}=t_a+p\Delta t$ with $p=\max\{n_0-1,2m_0-1\}$ before switching to the target step size $\Delta t$ for the remaining simulation until $t_{\max}$, cf. \cref{alg:ide_solver}. The finer step size is set according to~\eqref{eq:kh0} to $\Delta t_0 = \Delta t/k$ with $k=\left\lceil\Delta t^{-(n_0-1)/2}\right\rceil$. The results are depicted in~\cref{fig:adapted_timestep}.

In~\cref{fig:adapted_timestep} (top), we observe the second and third order convergence for Gregory orders 1 and 2 but also only third order convergence for Gregory rule 3. In the case of insufficient starting values, we have to increase the Gregory stepwise, thus we can apply Gregory rule 2 only from $t_2$ on, and Gregory rule 3 from $t_3$ on. This explains the lower convergence order for Gregory rule 3. For Gregory rule 2 we would expect only second order convergence, i.e., our observations are better than expected. We did not further investigate this.

In~\cref{fig:adapted_timestep}(bottom), we observe convergence order $n_0+1$ for Gregory rule $n_0$ until the error approaches machine precision. For $\Delta t=10^{-3}$, the error in the infected compartment is minimally larger than in the baseline scenario, cf.~\cref{fig:baseline} where the initial setting was not considered. This slight difference probably stems from the constants that were not considered when deriving the factor $k$ in~\eqref{eq:kh0}. By using the finer step size until $t_{switch}$, we can use Gregory rule from time point $t_{switch}$ on. Hence, at the first time point with respect to the target step size, we already obtain a sufficiently accurate approximation which is why we obtain the desired convergence with respect to the target step size.

\subsection{Impact of the included history on epidemic outcomes}\label{sec:impact_epidemic}

The advantage of IDE-based formulations over ODE-based ones clearly is the inclusion of past developments as these shape future dynamics. On one hand, including the history from the very first cases, i.e., $\widetilde{T}\approx\widehat{T}$ also means a substantially increased computational work load as $t-\widetilde{T}$ grows. On the other hand, setting $\widetilde{T}$ with increasing distance to $\widehat{T}$ will increase the model error as provided by~\eqref{eq:history_error_Sprime}; if we are not in the special case of exponential survival. In order to study the resulting model errors, we again generate ground truth trajectories. We consider two particular transition distributions, exponential and Erlang.

The general set up follows the baseline scenario as described in~\cref{sec:baseline}. We use a constant time step $\Delta t = 10^{-2}$ to approximate the continuous IDE solution. We set $\widehat{T}=-50$ and compare different values $\widetilde{T}\in\{-50,-15,0\}$ to see the impact of the model error through the assumption of initial infection age in the IDE simulations. In contrast to previous settings, we compute the different IDE solutions directly from $t_{0,\text{IDE}}=\widetilde{T}$ on.  We use the third Gregory quadrature rule to approximate the integrals and the fourth order FD scheme~\eqref{eq:FD4} to approximate the derivative $\widehat{S}'$. 

\subsubsection{Exponentially distributed sojourn times}

\begin{figure}[!ht]
	\centering
	\includegraphics[width=0.9\linewidth]{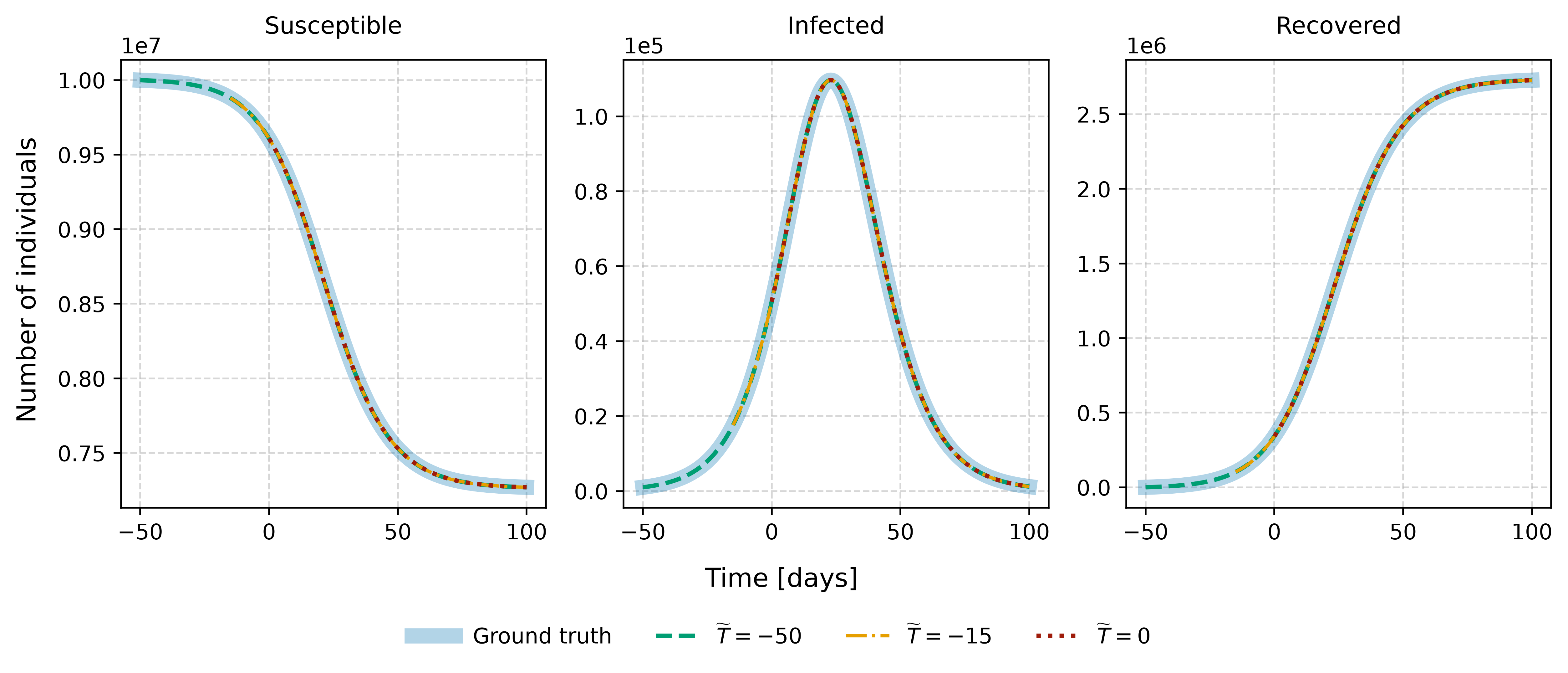}
	\caption{\textbf{Comparison of compartmental values for the ground truth and the three IDE simulations with different start values $\widetilde{T}$ using an exponentially distributed transition from infected to recovered state.} Numbers of individuals in susceptible, infected, and recovered state over the whole simulation period.}
	\label{fig:compare_compartments_exp}
	% Results from 2026-10-16
\end{figure}

\begin{figure}[!ht]
	\centering
	\includegraphics[width=0.9\linewidth]{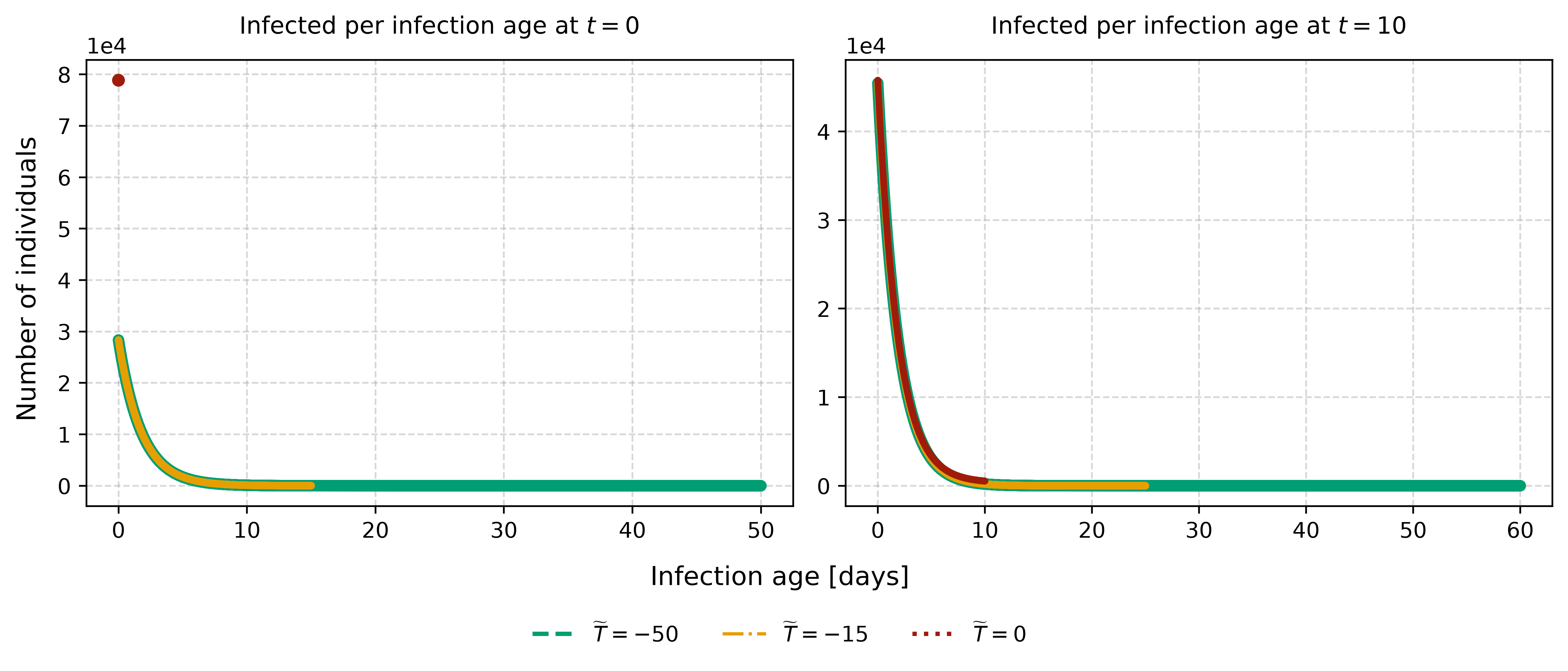}
	\caption{\textbf{Comparison of the infection age distribution for the three different start values $\widetilde{T}$ using an exponentially distributed transition from infected to recovered state.} \textbf{Left:} $t=0$, \textbf{right:} $t=10$.}
	\label{fig:compare_init_conditions_exp}
	% Results from 2026-10-16
\end{figure}

For all choices of $\widetilde{T}$, we see that the simulation results regarding the compartments coincide with the ground truth cf.~\cref{fig:compare_compartments_exp}. Regarding the infection age distribution, we have a high number of infected individuals with infection age zero at $t=0$ for $\widetilde{T}=0$ and a more diverse infection age distribution for the other two scenarios, cf.~\cref{fig:compare_init_conditions_exp} (left). This is due to assumption~\cref{eq:assump_initial_conditions}. At time $t=10$, we observe in~\cref{fig:compare_init_conditions_exp} (right) that the infection age distributions are much more aligned with only a small difference for an infection age of 10 as the remainder of the finite history assumption. In~\cref{lemma:ode,lemma:error}, we have shown that the IDE model is equivalent to a corresponding ODE model for any choice of $\widetilde{T}$ where the dynamics only depend on the current state of the system and not on the history. Our numerical simulations confirm these theoretical results. 

\subsubsection{Erlang-distributed sojourn times}
In contrast to the memoryless exponential distribution, the Erlang distribution is more realistically modeling compartmental transition~\cite{feng_epidemiological_2007}. To define a corresponding ground truth, we can use an ODE-based model employing the Linear Chain Trick (LCT model); cf., e.g.,~\cite{ploetzke_lct_2025} and the references therein.

In both, the LCT and the IDE, models, we consider Erlang distributed sojourn times with mean $T_I=4$. In the LCT model, we divide the infected compartment into six subcompartments. To be precise, the Erlang distribution is defined via the shape parameter $k=6$ and the rate parameter $\lambda = \frac{k}{T_I}$. At the earliest time point $\widehat{T}$, all infected individuals $\widehat{I}(\widehat{T})$ are assigned to the first subcompartment to match the assumption on infection age zero at $\widehat{T}$. Furthermore, we set $\phi=0.4$ and $t_\text{max}=80$. All other parameters follow~\cref{sec:baseline}.

\begin{figure}[!ht]
	\centering
	\includegraphics[width=0.9\linewidth]{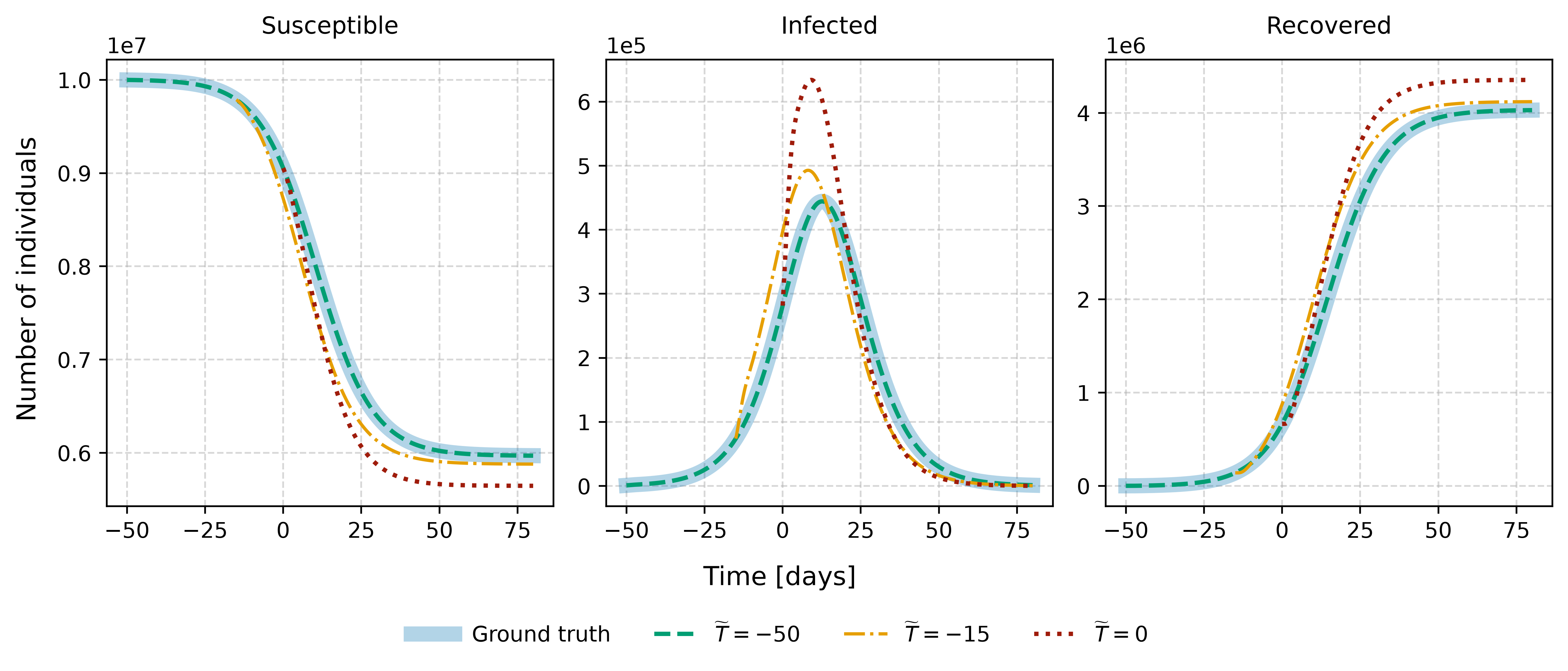} 
	\includegraphics[width=0.9\linewidth]{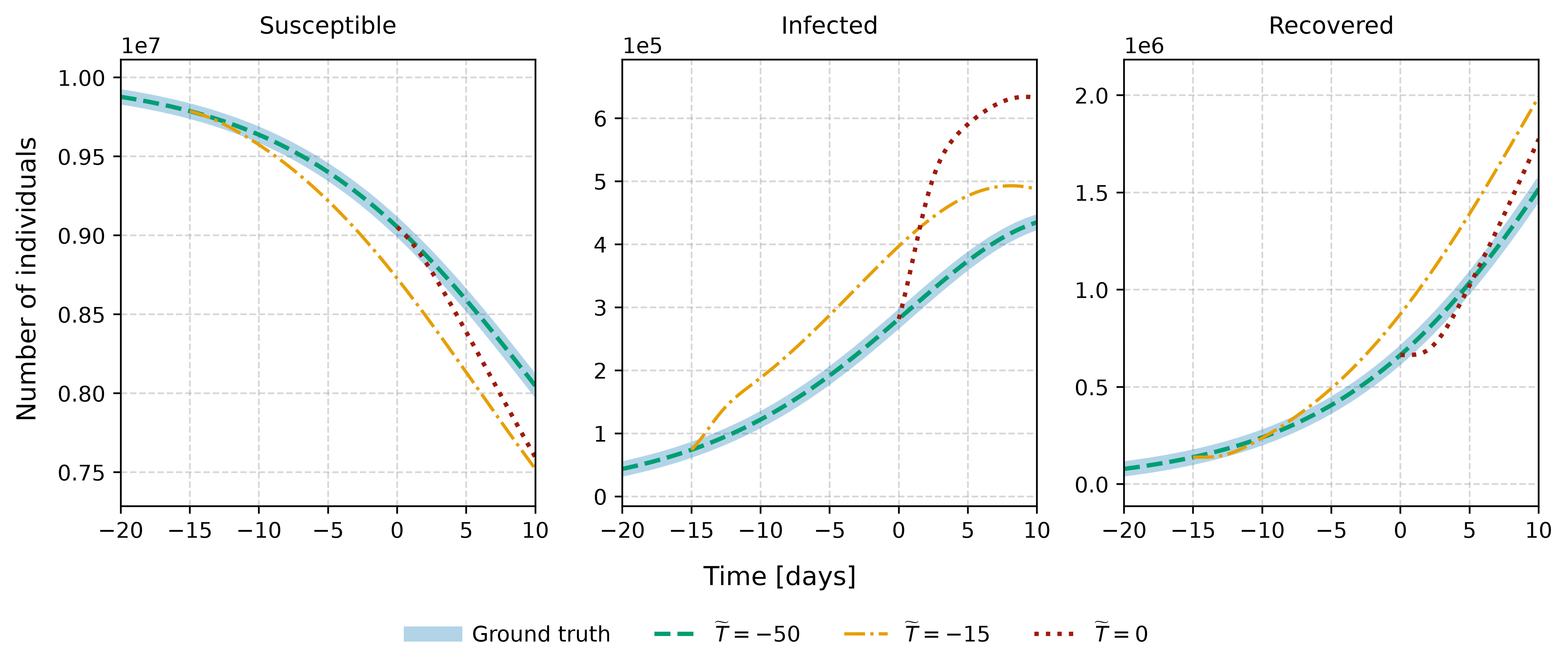} 
	\caption{\textbf{Comparison of compartmental values for the ground truth and the three IDE simulations with different start values $\widetilde{T}$ using an Erlang distributed transition from infected to recovered state.} \textit{Top row:} Numbers of individuals in susceptible, infected, and recovered state over the whole initialization and simulation period. 
		\textit{Bottom row:} Zoomed view of the top row from $t=-20$ until $t=10$.
	}
	\label{fig:compare_compartments_erlang}
	% Results from 2026-10-14
\end{figure}

In contrast to the memoryless exponential distribution, we see that the inclusion of historic data over differently sized initialization horizons can make a fundamental impact regarding peak size and timing; see~\cref{fig:compare_compartments_erlang}. Using $\widetilde{T}=\widehat{T}=-50$, the corresponding IDE simulation matches the ground truth very well, with only a difference of $1.79\cdot10^{-7}$ in the infected compartment at $t=10$. In contrast, when setting $\widetilde{T}=-15$, we obtain 52\,203.91 more infected individuals (+12.00~\% compared to the ground truth), and when setting $\widetilde{T}=0$, we obtain 197\,020.08 more infected individuals (+45.30~\% compared to the ground truth) at $t=10$, respectively.
Interestingly, we also observe that not only different peak sizes but also different final sizes are obtained for different values of $\widetilde{T}$ when using an Erlang distribution.

\begin{figure}[!ht]
	\centering
	\includegraphics[width=0.9\linewidth]{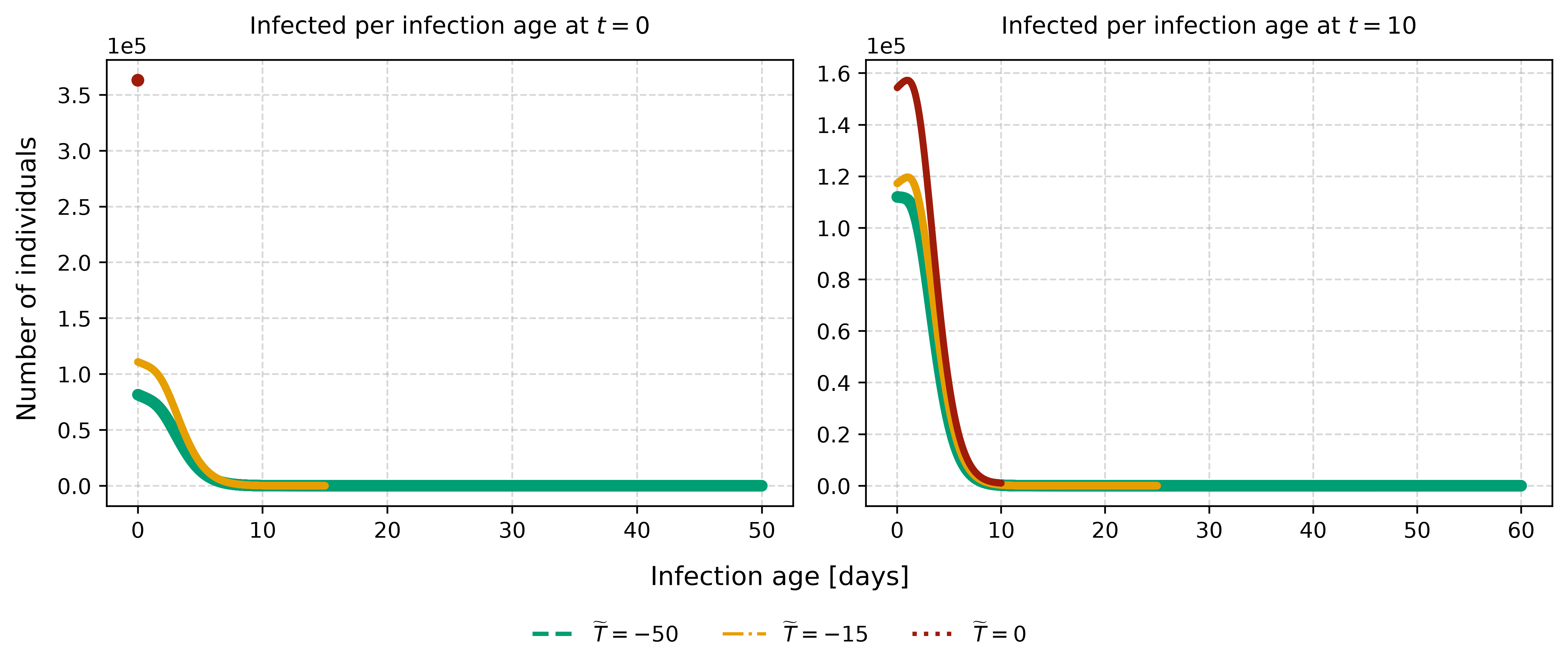}
	\caption{\textbf{Comparison of the infection age distribution for the three different start values $\widetilde{T}$ using an IDE model with Erlang distributed transition from infected to recovered state.} \textbf{Left:} $t=0$, \textbf{right:} $t=10$. }
	\label{fig:compare_init_conditions_erlang}
	% Results from 2026-10-16
\end{figure}

When using $\widetilde{T}>\widehat{T}$, we observe more infected individuals in the rise of the epidemic than in the ground truth, cf.~\cref{fig:compare_compartments_erlang} (bottom). By making the zero infection age assumption at $\widetilde{T}>\widehat{T}$, we are overestimating the number of individuals that have a small infection age for which we have larger values of $\gamma_I^R$; cf.~\cref{fig:compare_init_conditions_erlang}. Hence, we are overestimating the contribution to the force of infection these individuals make, thus leading to more infected individuals than in the ground truth. 

Moreover, we see that the number of recovered individuals stays nearly constant for the first days after $\widetilde{T}$ when setting $\widetilde{T}>\widehat{T}$, cf.~\cref{fig:compare_compartments_erlang} (bottom). All infected individuals are assigned an infection age of zero at $\widetilde{T}$, cf.~\eqref{eq:assump_initial_conditions}. Since $\gamma_I^R(0)=1$, only a small fraction of these individuals is recovering in the first days. In the ground truth, the infection age distribution is more diverse at the respective time points $t\in\{-15,0\}$ so that more individuals have a larger infection age thus leading to more individuals that are recovering in the first days after $t\in\{-15,0\}$. For this reason, we observe a delay in the rise of the recovered compartment in the scenarios with $\widetilde{T}>\widehat{T}$ compared to the ground truth.

\section{Conclusion}\label{sec:conclusion}

In this paper, we considered different aspects of the numerical solution of IDE-based models for infectious disease dynamics. For similar IDE-based formulations selected theoretical properties of the model~\cite{brauerAgeinfectionFinalSize2008,brauerDynamicsAgeinfectionCholera2013} or the novel numerical schemes~\cite{messina_positive_2022,pezzellaInfectiousDisease2024} were considered. However, when it came to the application of IDE-based models for ongoing epidemic dynamics, several questions remained unanswered. We generalized the model formulation in the number of susceptible individuals by allowing contact change points in the simulation. Importantly, we also considered the developments and convergence rates in further disease compartments such as infected and recovered. 

The numerical extension for additional compartments requires additional components, e.g., particular finite differences to approximate the derivative in the number of susceptible individuals. Computing additional compartments is of particular importance for practical applications, e.g., to also obtain the total number of ongoing infections or the number of new hospitalizations in more advanced models~\cite{wendler_nonstandard_2026}. While the convergence in the number of susceptibles only depends on the Gregory rule, we showed that similar, higher order FD schemes are needed to obtain the intended convergence order in the number of infected or recovered individuals.

The characterization of the time-dependent contact rate as a step-wise constant function with $\mathcal{C}^k$ smoothing between the jumps allows for setting and calibrating contacts on a finite number of discretization points. We compared a previously implemented cosine smoother function of the MEmilio framework~\cite{bicker2026memilio} against polynomial smoothstep functions that yield up to fourth order convergence for settings with intervention-induced contact change points. An extension to even higher order methods is straightforward as it only requires replacing the individual components of our scheme (Gregory rule, FD, and smoothstep function) by functions with higher order error terms.

When comparing to a ground truth obtained with an ODE model, we showed for the exponential setting that the convergence order of the higher order numerical scheme deteriorates if we do not have sufficiently accurate approximations of past dynamics. We showed that the requirement of suitable initial values can be resolved by using a finer step size on a small initial interval before switching to a larger time step. 

We considered IDE models that were started not only at the start of the epidemic but also at later times in the outbreak. This required making the assumption of zero infection age at the simulation start. 
As expected, we showed that starting the simulation at different time points does not change the course of the epidemic if an exponential transition distribution is used. This is in contrast to more realistic transition distributions such as Erlang for which we have shown substantial differences in epidemic outcomes when starting the IDE simulation at different times in the epidemic. 

The IDE-based SIR model was implemented in C++ in the high performance Modular EpideMIcs simuLatIOn software MEmilio. While we demonstrated how the solver's individual components yield schemes with up to fourth order convergence, we also demonstrated that less costly components such as, e.g., lower order FD are sufficient if lower order Gregory rules are used. As the MEmilio framework shall enable users to flexibly choose between models and numerical solution schemes, the next software engineering tasks will contain the adaptive and compile-time combination of the individual solver components to obtain a run-time to convergence order optimized numerical solution scheme. In here, comparisons will have to be made to the already existing lower order nonstandard scheme for which a cost-efficient update rule was established~\cite{wendler_nonstandard_2026}. We furthermore expect our scheme to be generalizable to, e.g., SECIR-type models with presymptomatic or hospitalized individuals.

\section*{Author Contributions}

\noindent\textbf{Conceptualization:} Anna Wendler, Martin J. Kühn\\
\noindent\textbf{Data Curation:} Anna Wendler\\
\noindent\textbf{Formal Analysis:} Anna Wendler, Jan Hasenauer, Martin J. Kühn\\
\noindent\textbf{Funding Acquisition:} Jan Hasenauer, Martin J. Kühn\\
\noindent\textbf{Investigation:} Anna Wendler, Martin J. Kühn\\
\noindent\textbf{Methodology:} Anna Wendler, Martin J. Kühn\\
\noindent\textbf{Project Administration:} Martin J. Kühn\\
\noindent\textbf{Resources:} Martin J. Kühn\\
\noindent\textbf{Software:} Anna Wendler\\
\noindent\textbf{Supervision:} Jan Hasenauer, Martin J. Kühn \\
\noindent\textbf{Validation:} Anna Wendler, Jan Hasenauer, Martin J. Kühn \\
\noindent\textbf{Visualization:} Anna Wendler, Martin J. Kühn\\
\noindent\textbf{Writing – Original Draft:} Anna Wendler, Martin J. Kühn\\
\noindent\textbf{Writing – Review \& Editing:} Anna Wendler, Jan Hasenauer, Martin J. Kühn

%\newpage
%\pagenumbering{roman}
\appendix

\section*{Acknowledgements}
This work was supported by the Initiative and Networking Fund of the Helmholtz Association (grant agreement number KA1-Co-08, Project LOKI-Pandemics). It was furthermore supported by the German Federal Ministry of Education and Research and the German Federal Ministry of Research, Technology and Space under grant agreement 031L0324E (Project EPISERVE). EPISERVE is part of the Modeling Network for Severe Infectious Diseases (MONID). Additionally, this work was supported by the European Union via the ERC grant INTEGRATE, grant agreement number 101126146, and under Germany’s Excellence Strategy by the Deutsche Forschungsgemeinschaft (DFG, German Research Foundation) (EXC 2047—390685813, EXC 2151—390873048, and 524747443), the University of Bonn via the Schlegel Professorship of J.H. 

\section*{Competing interests}
The authors declare that they have no known competing financial interests or personal relationships that could have appeared to influence the work reported in this paper.

\section*{Data availability}
All simulations were implemented and performed using the MEmilio framework \cite{bicker2026memilio}. The specific simulation scenarios and plotting scripts are accessible at \url{https://github.com/SciCompMod/memilio-simulations/tree/main/2026_Wendler_et_al_Higher_order_IDE}.

\section{Gregory quadrature weights}\label{app:gregory_weights}

We consider an integral approximation of the form  
\begin{align*}
	\int_{\widetilde{T}}^{t_n} f(x)\, \d x \approx \Delta t \sum_{j=a}^n \omega_{n-a,j-a}^{n_0} \, f(t_j),\quad n\geq a+n_0,
\end{align*}
where $t_k=k\,\Delta t$ for $k\in\mathbb{Z}$ provides a uniform discretization of $\mathbb{R}$, and we assume $t_{a}=\widetilde{T}$, $a\in\mathbb{Z}_{\leq 0}$. 

The quadrature weights for Gregory rule $n_0$ can be described via the matrices $\Omega^{n_0}=\Big(\omega_{i,j}^{n_0}\Big)_{i,j}$ with $i\in\mathbb{N}_{\geq n_0}$, $j\in\{0,\dots,i\}$. For each $\Omega^{n_0}$, row $i$ contains $i+1$ weights that are correspondingly used for an approximation based on $i+1$ discrete points, i.e., an approximation through a sum with $i+1$ summands. The column index $j$ refers to the weight of the corresponding summand within this sum. Note that for Gregory rule $n_0$, the weights for the first and last $n_0$ discrete points are corrected whereas the remaining inner points have a weight of one and that the correction is symmetric. 
For the derivation of the specific weights $\omega_{i,j}^{n_0}$, we refer to~\cite{brunnerNumericalSolutionVolterra1986}. Examples for $n_0\in\{1,2,3\}$, which are used within this paper, are given below. 

\begin{itemize}
	\item Gregory rule 1:
	\begin{align}
		\Omega^1 = \frac{1}{2}\,
		\begin{pmatrix}
			1 & 1 \\
			1 & 2 & 1 \\
			\vdots & \vdots & \ddots & \ddots \\
			1 & 2 & \cdots & 2 & 1 \\
		\end{pmatrix}
	\end{align}
	\item Gregory rule 2:
	\begin{align}
		\Omega^2 = \frac{1}{12}\,
		\begin{pmatrix}
			5 & 14 & 5 \\
			5 & 13 & 13 & 5 \\
			5 & 13 & 12 & 13 & 5 \\
			% 5 & 13 & 12 & 12 & 13 & 5 \\
			\vdots & \vdots & \vdots& \ddots & \ddots & \ddots \\
			5 & 13 & 12 & \cdots & 12 & 13 & 5
		\end{pmatrix}
	\end{align}
	\item Gregory rule 3
	\begin{align}
		\Omega^3 = \frac{1}{24}\,
		\begin{pmatrix}
			9 & 27 & 27 & 9 \\
			9 & 28 & 22 & 28 & 9 \\
			9 & 28 & 23 & 23 & 28 & 9 \\
			9 & 28 & 23 & 24 & 23 & 28 & 9 \\
			\vdots & \vdots & \vdots& \vdots &  \ddots & \ddots & \ddots & \ddots \\
			9 & 28 & 23 & 24 & \cdots & 24 & 23 & 28 & 9
		\end{pmatrix}
	\end{align}
\end{itemize} 

\section{Polynomial smoothstep functions}\label{app:smoothstep_def}

In~\cref{fig:runtime_smoothstep}, we consider not only the cosine smoother cosine and $\mathcal{C}^4$ polynomial smoothstep functions, but also $\mathcal{C}^k, k\in\{1,2,3\}$ smoothstep functions. These functions are defined below where $\phi_0$, $t_D$ and $\delta$ are as in~\cref{sec:continuous_model}.

\begin{align}
	\widetilde{\phi}^{P_1}(t) &= \phi_0 + ((1-r)\,\phi_0 - \phi_0) \left(3\left(\frac{t - t_D+ \delta}{\delta}\right)^2 -2\left(\frac{t - t_D+ \delta}{\delta}\right)^3\right),\\
	\widetilde{\phi}^{P_2}(t) &= \phi_0 + ((1-r)\,\phi_0 - \phi_0) \left(10\left(\frac{t - t_D+ \delta}{\delta}\right)^3 -15\left(\frac{t - t_D+ \delta}{\delta}\right)^4 + 6 \left(\frac{t - t_D+ \delta}{\delta}\right)^5\right),\\
	\widetilde{\phi}^{P_3}(t) &= \phi_0 + ((1-r)\,\phi_0 - \phi_0) \left(35\left(\frac{t - t_D+ \delta}{\delta}\right)^4 -84\left(\frac{t - t_D+ \delta}{\delta}\right)^5 + 70 \left(\frac{t - t_D+ \delta}{\delta}\right)^6 -20\left(\frac{t - t_D+ \delta}{\delta}\right)^7\right).
\end{align}

\bibliographystyle{elsarticle-num}
\bibliography{literature}

\end{document}